\documentclass[a4paper]{amsart} 

\usepackage{enumerate}
\usepackage{amssymb}
\usepackage{bbold}
\usepackage[all]{xy}
\SelectTips{cm}{}

\usepackage[colorlinks]{hyperref}
\calclayout 
\makeatletter
\def\serieslogo@{} 
\def\@setcopyright{} 

\makeatother

\title[Local cohomology and support]{Local cohomology and support for \\ triangulated categories}  
\thanks{Version from 11th February 2008.}  
\thanks{D.B.\ was supported by a Senior Scientist Prize from the Humboldt Foundation.  S.I.\ was supported by NSF grant, DMS 0602498} 
\author{Dave Benson}
\address{Dave Benson \\
Department of Mathematical Sciences\\
University of Aberdeen\\
Meston Building\\
King's College\\
Aberdeen AB24 3UE\\
Scotland U.K.} 

\author{Srikanth B. Iyengar}
\address{Srikanth B. Iyengar\\ Department of Mathematics\\
University of Nebraska\\ Lincoln NE 68588\\ U.S.A.}
\email{iyengar@math.unl.edu}

\author{Henning Krause}
\address{Henning Krause\\ Institut f\"ur Mathematik\\
Universit\"at Paderborn\\ 33095 Paderborn\\ Germany.}
\email{hkrause@math.upb.de}

\dedicatory{To Lucho Avramov, on his 60th birthday.}

\theoremstyle{plain}

\newtheorem{theorem}{Theorem}[section]

\newtheorem{proposition}[theorem]{Proposition}

\newtheorem{lemma}[theorem]{Lemma}

\newtheorem{corollary}[theorem]{Corollary}

\newtheorem{itheorem}{Theorem}

\theoremstyle{definition}

\newtheorem{definition}[theorem]{Definition}

\newtheorem{example}[theorem]{Example}

\newtheorem{notation}[theorem]{Notation}
\newtheorem*{ack}{Acknowledgments}

\theoremstyle{remark}

\newtheorem{remark}[theorem]{Remark}
\newtheorem*{Claim}{Claim}

\numberwithin{equation}{theorem}

 \newcommand{\col}{\colon}

\newcommand{\Loc}{\operatorname{Loc}\nolimits}

\renewcommand{\dim}{\operatorname{dim}\nolimits}
\newcommand{\End}{\operatorname{End}\nolimits}
\newcommand{\Ext}{\operatorname{Ext}\nolimits}
\newcommand{\Tor}{\operatorname{Tor}\nolimits}
\newcommand{\funct}{\operatorname{\mathcal{H}\!\!\;\mathit{om}}\nolimits}
\newcommand{\Hom}{\operatorname{Hom}\nolimits}

\newcommand{\id}{\operatorname{id}\nolimits}
\newcommand{\Id}{\operatorname{Id}\nolimits}

\renewcommand{\Im}{\operatorname{Im}\nolimits}
\newcommand{\Ker}{\operatorname{Ker}\nolimits}
\newcommand{\cone}{\operatorname{cone}\nolimits}

\newcommand{\ann}{\operatorname{ann}\nolimits}
\newcommand{\spec}{\operatorname{Spec}\nolimits}
\newcommand{\supp}{\operatorname{supp}\nolimits}
\newcommand{\Min}{\operatorname{min}\nolimits}

\renewcommand{\mod}{\operatorname{\mathsf{mod}}\nolimits}
\newcommand{\Mod}{\operatorname{\mathsf{Mod}}\nolimits}
\newcommand{\StMod}{\operatorname{\mathsf{StMod}}\nolimits}
\newcommand{\stmod}{\operatorname{\mathsf{stmod}}\nolimits}
\newcommand{\Inj}{\operatorname{\mathsf{Inj}}\nolimits}
\newcommand{\height}{\operatorname{ht}\nolimits}
\newcommand{\cl}{\operatorname{cl}\nolimits}

\newcommand{\hH}{H}
\newcommand{\hh}[1]{{H(#1)}}
\newcommand{\HH}{\mathrm{HH}}
\newcommand{\coh}[3][*]{{H^{#1}_{#2}(#3)}}
\newcommand{\ac}{\mathrm{ac}}

\newcommand{\inc}{\mathrm{inc}}
\newcommand{\op}{\mathrm{op}}

\newcommand{\lto}{\longrightarrow}
\newcommand{\xra}{\xrightarrow}

\newcommand{\comp}{\mathop{\circ}}

\newcommand{\un}{\un}

\newcommand{\wt}{\widetilde} 

\def\Si{{\Sigma\,}}
\def\si{\sigma}

\def\mcU{\mathcal{U}}
\def\mcV{\mathcal{V}}
\def\mcW{\mathcal{W}}
\def\mcY{\mathcal{Y}}
\def\mcZ{\mathcal{Z}}
\def\bbz{\mathbb Z}
\def\bbi{\mathbb 1}

\def\sfe{\mathsf e}

\def\sfi{\mathsf i}

\def\sfq{\mathsf q}

\def\sfA{\mathsf A}
\def\sfC{\mathsf C}
\def\sfD{\mathsf D}

\def\sfG{\mathsf G}
\def\sfM{\mathsf M}
\def\sfS{\mathsf S}
\def\sfT{\mathsf T}
\def\sfK{\mathsf K}

\def\bfr{\mathbf{r}}

\newcommand{\bsq}{\boldsymbol{q}}
\newcommand{\eps}{\varepsilon}
\newcommand{\fa}{\mathfrak{a}}
\newcommand{\fm}{\mathfrak{m}}
\newcommand{\fp}{\mathfrak{p}}
\newcommand{\fq}{\mathfrak{q}}

\newcommand{\lch}[1]{{{\mathbf R}F}_{#1}}

\newcommand{\gam}{\varGamma}
\newcommand{\kos}[2]{{#1/\!\!/#2}}
\newcommand{\bloc}{{L}}
\newcommand{\cent}{{Z}}

\begin{document}

\begin{abstract}
We propose a new method for defining a notion of support for objects
in any compactly generated triangulated category admitting small
coproducts.  This approach is based on a construction of local
cohomology functors on triangulated categories, with respect to a
central ring of operators.  Suitably specialized one recovers, for
example, the theory for commutative noetherian rings due to Foxby and
Neeman, the theory of Avramov and Buchweitz for complete intersection
local rings, and varieties for representations of finite groups
according to Benson, Carlson, and Rickard. We give explicit examples
of objects whose triangulated support and cohomological support
differ. In the case of group representations, this allows us to
correct and establish a conjecture of Benson.
\end{abstract} 
\maketitle

\setcounter{tocdepth}{1}
\tableofcontents

\section{Introduction}

\hfill \begin{minipage}[t]{6.5cm} \small Herr K.\ sagte einmal: ``Der Denkende ben\"utzt
  kein Licht zuviel, kein St\"uck Brot zuviel, keinen Gedanken zuviel.''  \smallskip
  
\noindent \footnotesize {\sc Bertolt Brecht}, Geschichten von Herrn Keuner
\end{minipage}

\normalsize

\vspace*{1em}

The notion of support is a fundamental concept which provides a geometric approach for
studying various algebraic structures. The prototype for this has been Quillen's~\cite{Qu}
description of the algebraic variety corresponding to the cohomology ring of a
finite group, based on which Carlson~\cite{Ca:1983} introduced \emph{support} varieties for
modular representations. This has made it possible to apply methods of algebraic geometry
to obtain representation theoretic information.  Their work has inspired the development
of analogous theories in various contexts, notably modules over commutative complete
intersection rings, and over cocommutative Hopf algebras.

In this article we propose a new method for defining a notion of support for objects in
any compactly generated \emph{triangulated category} admitting small coproducts.  The
foundation of our approach is a construction of \emph{local cohomology} functors on
triangulated categories, with respect to a central ring of operators; this is inspired by
work of Grothendieck~\cite{Gr}.  Suitably specialized our approach recovers, for example,
the support theory of Foxby~\cite{Fo} and Neeman~\cite{Ne:dr} for commutative noetherian
rings, the theory of Avramov and Buchweitz for complete intersection local
rings~\cite{Av:vpd,AB}, and varieties for representations of finite groups, according to
Benson,Carlson, and Rickard~\cite{BCR:1996a}. It is surprising how little is needed to
develop a satisfactory theory of support. To explain this, let us sketch the main results
of this paper.

Let $\sfT$ be a triangulated category that admits small coproducts and is compactly
generated. In the introduction, for ease of exposition, we assume $\sfT$ is generated by a
single compact object $C_0$. Let $Z(\sfT)$ denote the graded center of $\sfT$. The notion
of support presented here depends on the choice of a graded-com\-mu\-ta\-tive noetherian
ring $R$ and a homomorphism of rings
\[
R\lto Z(\sfT)\,.\]
  We may view $R$ as a ring of cohomology operators on $\sfT$. For each object $X$ in $\sfT$
its cohomology
\[
H^*(X)=\Hom_\sfT^*(C_0,X)=\coprod_{n\in\mathbb Z}\Hom_\sfT(C_0,\Si^nX)\]
  has a structure of a graded module over $Z(\sfT)$ and hence over $R$. We let $\spec R$
denote the set of graded prime ideals of $R$.  The \emph{specialization closure} of a
subset $\mcU\subseteq\spec R$ is the subset
\[
\cl\mcU=\{\fp\in\spec R\mid\text{there exists $\fq\in\mcU$ with $\fq\subseteq \fp$}\}\,.\]
  This is the smallest specialization closed subset containing $\mcU$.

One of the main results of this work is an axiomatic characterization of support:

\begin{itheorem}\label{intro:main}
There exists a unique assignment sending each object $X$ in $\sfT$ to a subset $\supp_{R}
X$ of $\spec R$ such that the following properties hold:
\begin{enumerate}[{\quad\rm(1)}]
\item \emph{Cohomology:} For each object $X$ in $\sfT$ one has
\[
\cl(\supp_{R} X) = \cl(\supp_R H^*(X))\,.\]
  \item \emph{Orthogonality:} For objects $X$ and $Y$ in $\sfT$, one has that
\[
\cl(\supp_{R} X)\cap\supp_{R} Y=\varnothing \quad\text{implies}\quad \Hom_\sfT(X,Y)=0\,.\]
  \item \emph{Exactness:} For every exact triangle $W\to X\to Y\to$ in $\sfT$, one has
\[
\supp_{R} X\subseteq \supp_{R} W \cup \supp_{R} Y\,.\]
  \item \emph{Separation:} For any specialization closed subset $\mcV$ of $\spec R$ and any
  object $X$ in $\sfT$, there exists an exact triangle $X'\to X\to X''\to$ in $\sfT$ such
  that
\[
\supp_{R} X'\subseteq\mcV \quad\text{and} \quad\supp_{R}X'' \subseteq\spec R\setminus \mcV\,.\]
  \end{enumerate}
\end{itheorem}

Implicit in (1) is a comparison of the \emph{triangulated support} $\supp_{R} X$ and the
\emph{cohomological support} $\supp_R H^*(X)$. This was part of the initial motivation for
this work. We prove also that if the cohomology $H^*(X)$ is finitely generated as a module
over $R$, then $\supp_{R} X=\supp_R H^*(X)$. Without such finiteness assumption however,
triangulated and cohomological support can differ; see Sections~\ref{se:commutative} and
\ref{se:groups} for explicit examples.

It is thus interesting that the triangulated support of an object $X$ can be yet detected
by cohomology. Only, one has to compute cohomology with respect to \emph{each compact
  object}.  This is made precise in the next result, where, for a graded $R$-module $M$,
we write $\Min_R M$ for the set of minimal primes in its support.

\begin{itheorem}\label{intro:support}
For each object $X$ in $\sfT$, one has an equality:
\[
\supp_{R} X = \bigcup_{C\text{ compact}} \Min_R \Hom_\sfT^*(C,X) \,.\]
  In particular, $\supp_{R} X=\varnothing$ if and only if $X=0$.
\end{itheorem}

Beyond proving Theorems 1 and 2, we develop systematically a theory of supports in order
to make it a viable tool.  For example, in Section~\ref{se:connectedness}, we establish
the following result of Krull-Remak-Schmidt type.

\begin{itheorem}\label{intro:KRS}
Each object $X$ in $\sfT$ admits a unique decomposition $X=\coprod_{i\in I} X_i$ with
$X_i\ne 0$ such that the subsets $\cl(\supp_{R} X_i)$ are connected and pairwise disjoint.
\end{itheorem}

Here is a direct corollary: If $X$ is an indecomposable object in $\sfT$, then $\supp_{R}
X$ is a connected subset of $\spec R$. This generalizes and unifies various connectedness
results in the literature, starting with a celebrated theorem of Carlson, which states
that the variety of an indecomposable group representation is connected \cite{Ca}.

As stated before, the basis for this work is a construction of local
cohomology functors on $\sfT$. Given a specialization closed subset
$\mcV$ of $\spec R$, we establish the existence of (co)localization
functors $\gam_{\mcV}$ and $\bloc_{\mcV}$ on $\sfT$, such that for
each $X$ in $\sfT$ there is a natural exact triangle
\[
\gam_{\mcV}X\lto X\lto \bloc_{\mcV}X\lto \]
  in $\sfT$. We view $\gam_{\mcV}$ as the local cohomology functor with respect to $\mcV$.
One justification for this is the following result:
\[
\supp_{R}\gam_{\mcV}X = \mcV \cap \supp_{R}X \quad\text{and}\quad
\supp_{R}\bloc_{\mcV}X = \big(\spec R\setminus \mcV\big) \cap \supp_{R}X \,.\]
  
A major focus of this work are properties of the functors $\gam_{\mcV}$ and $\bloc_{\mcV}$
for a general triangulated category $\sfT$; the results on support are derived from
them. These occupy Sections 4--7 in this article; the first three prepare the ground for
them, and for later sections.  The remaining sections are devoted to various specific
contexts, and are intended to demonstrate the range and applicability of the methods
introduced here.  We stress that hitherto many of the results established were known only
in special cases; Theorem~\ref{intro:KRS} is such an example. Others, for instance,
Theorems~\ref{intro:main} and \ref{intro:support}, are new in all contexts relevant to
this work.

In Section~\ref{se:tensor} we consider the case where the triangulated category $\sfT$
admits a symmetric tensor product. The notion of support then obtained is shown to
coincide with the one introduced by Hovey, Palmieri, and Strickland~\cite{HPS}.

Section~\ref{se:commutative} is devoted to the case where $\sfT$ is the derived category
of a commutative noetherian ring $A$, and $R=A\to Z(\sfT)$ is the canonical morphism.  We
prove that for each specialization closed subset $\mcV$ of $\spec A$ and complex $X$ of
$A$-modules the cohomology of $\gam_{\mcV}X$ is classical local cohomology, introduced by
Grothendieck \cite{Gr}.

The case of modules for finite groups is studied in Section~\ref{se:groups}, where we
prove that support as defined here coincides with one of Benson, Carlson and
Rickard~\cite{BCR:1996a}. Even though this case has been studied extensively in the
literature, our work does provide interesting new information. For instance, using
Theorem~\ref{intro:support}, we describe an explicit way of computing the support of a
module in terms of its cohomological supports.  This, in spirit, settles Conjecture 10.7.1
of \cite{Benson:2004b} that the support of a module equals the cohomological support; we
provide an example that shows that the conjecture itself is false.

The final Section~\ref{se:ci} is devoted to complete intersection local rings. We recover
the theory of Avramov and Buchweitz for support varieties of finitely generated modules
\cite{Av:vpd,AB}. A salient feature of our approach is that it gives a theory of local
cohomology with respect to rings of cohomology operators.

This article has influenced some of our subsequent work on this topic:
in \cite{AI}, Avramov and Iyengar address the problem of realizing
modules over arbitrary associative rings with prescribed cohomological
support; in \cite{Kr:thick}, Krause studies the classification of
thick subcategories of modules over commutative noetherian
rings. Lastly, the techniques introduced here play a pivotal role in
our recent work on a classification theorem for the localizing
subcategories of the stable module category of a finite group; see
\cite{BIK:2008}.

\section{Support for modules}\label{se:modules}

In this section $R$ denotes a $\bbz$-graded-commutative noetherian ring. Thus we have
$x\cdot y=(-1)^{|x||y|}y\cdot x$ for each pair of homogeneous elements $x,y$ in $R$.

Let $M$ and $N$ be graded $R$-modules. For each integer $n$, we write $M[n]$ for the
graded module with $M[n]^i = M^{i+n}$. We write $\Hom^*_R(M,N)$ for the graded
homomorphisms between $M$ and $N$:
\[ 
\Hom^n_R(M,N)= \Hom_R(M,N[n])\,. 
\]
  
The degree zero component is usually abbreviated to $\Hom_R(M,N)$.  Since $R$ is
graded-commutative, $\Hom^*_R(M,N)$ is a graded $R$-module in an obvious way.  Henceforth,
unless otherwise specified, when we talk about modules, homomorphisms, and tensor
products, it is usually implicit that they are graded.  We write $\Mod R$ for the category
of graded $R$-modules.

\subsection*{Spectrum}
Let $\spec R$ denote the set of graded prime ideals of $R$.  Given a homogeneous ideal
$\fa$ in $R$, we set:
\[
\mcV(\fa) = \{\fp\in\spec R\mid \fp\supseteq \fa\}\,.\]
  
Such subsets of $\spec R$ are the closed sets in the \emph{Zariski topology} on $\spec R$.
Let $\mcU$ be a subset of $\spec R$. The specialization closure of $\mcU$ is the set
\[
\cl\mcU=\{\fp\in\spec R\mid\text{there exists $\fq\in\mcU$ with $\fq\subseteq \fp$}\}\,.\]
  The subset $\mcU$ is \emph{specialization closed} if $\cl\mcU=\mcU$. Evidently,
specialization closed subsets are precisely the unions of Zariski closed subsets of $\spec R$.

Let $\fp$ be a prime ideal.  We write $R_\fp$ for the homogeneous localization of $R$ with
respect to $\fp$; it is a graded local ring in the sense of Bruns and Herzog~\cite[(1.5.13)]{BH}, with maximal
ideal $\fp R_{\fp}$. The graded field $R_\fp/\fp R_\fp$ is denoted $k(\fp)$.
As usual, we write $E(R/\fp)$ for the injective envelope of the $R$ module $R/\fp$.  
Given an $R$-module $M$ we let $M_\fp$ denote the localization of $M$ at $\fp$.

\subsection*{Injective modules}
The classification of injective modules over commutative noetherian rings carries over to
the graded case with little change: Over a graded-commutative noetherian ring $R$ an
arbitrary direct sum of injective modules is injective; every injective module decomposes
essentially uniquely as a direct sum of injective indecomposables. Moreover, for each
prime ideal $\fp$ in $R$, the injective hull $E(R/\fp)$ of the quotient $R/\fp$ is
indecomposable, and each injective indecomposable is isomorphic to a shifted copy
$E(R/\fp)[n]$ for a unique prime $\fp$ and some not necessarily unique $n\in \mathbb
Z$. For details, see \cite[(3.6.3)]{BH}.

\subsection*{Torsion modules and local modules}
Let $R$ be a graded-commutative noetherian ring, $\fp$ a prime ideal in $R$, and let $M$
be an $R$-module. The module $M$ is said to be \emph{$\fp$-torsion} if each element of $M$
is annihilated by a power of $\fp$; equivalently:
\[
M =\{m\in M\mid \text{there exists an integer $n\geq 0$ such that $\fp^n\cdot m=0$.}\}\]
  The module $M$ is \emph{$\fp$-local} if the natural map $M\to M_\fp$ is bijective.

For example, $R/\fp$ is $\fp$-torsion, but it is $\fp$-local only if $\fp$ is a maximal
ideal, while $R_\fp$ is $\fp$-local, but it is $\fp$-torsion only if $\fp$ is a minimal
prime ideal.  The $R$-module $E(R/\fp)$, the injective hull of $R/\fp$, is both
$\fp$-torsion and $\fp$-local. The following lemma is easy to prove.

\begin{lemma} 
\label{le:inj-local}
\pushQED{\qed}%
Let $\fp$ be a prime ideal in $R$. For each prime ideal $\fq$ in $R$ one has
\[
E(R/\fp)_\fq = \begin{cases} 
E(R/\fp)  & \text{if $\fq\in\mcV(\fp)$,}\\
0 & \text{otherwise.}
\end{cases} \qedhere
\]
  \end{lemma}

\subsection*{Injective resolutions}
Each $R$-module $M$ admits a \emph{minimal injective resolution}, and such a resolution is
unique, up to isomorphism of complexes of $R$-modules.  We say that $\fp$ occurs in a
minimal injective resolution $I$ of $M$, if for some pair of integers $i,n\in\mathbb Z$,
the module $I^i$ has a direct summand isomorphic to $E(R/\fp)[n]$. The \emph{support} of
$M$ is the set
\[
\supp_R M=\left\{\fp\in\spec R \left|
\begin{gathered}
\text{$\fp$ occurs in a minimal} \\
\text{ injective resolution of $M$}
\end{gathered}
\right.\right\}\]
  In the literature, $\supp_RM$ is sometimes referred to as the cohomological support, or
the small support, of $M$, in order to distinguish it from the usual support, which is the
subset $\{\fp\in \spec R\mid M_\fp\ne 0\}$, sometimes denoted $\mathrm{Supp}_RM$; see
\cite{Fo} and also Lemma~\ref{le:supp-ann}.  The cohomological support has other
descriptions; see Section~\ref{se:commutative}.

We now recollect some properties of supports. In what follows, the annihilator of an
$R$-module $M$ is denoted $\ann_RM$.

\begin{lemma}\label{le:supp-ann}
The following statements hold for each $R$-module $M$.
\begin{enumerate}[{\quad\rm(1)}]
\item One has inclusions
\[
\supp_RM\subseteq \cl(\supp_RM) =\{\fp\in \spec R\mid M_\fp\ne 0\}
     \subseteq\mcV(\ann_RM)\,,\]
  and equalities hold when $M$ is finitely generated.
\item
For each $\fp\in\spec R$, one has an equality
\[
\supp_{R}(M_\fp)=\supp_RM\cap \{\fq\in\spec R\mid \fq\subseteq \fp\}\,.\]
  \end{enumerate}
\end{lemma}

\begin{proof}
  If $I$ is a minimal injective resolution of $M$ over $R$, then $I_\fp$ is a minimal
  injective resolution of $M_\fp$, see \cite[\S18]{Mat}.  In view of
  Lemma~\ref{le:inj-local}, this implies (2) and the equality in (1).  The inclusions in
  (1) are obvious. 
  
It remains to verify that when $M$ is finitely generated $\mcV(\ann_RM)\subseteq \supp_RM$ holds.
For this, it suffices to prove that for any prime ideal $\fp\supseteq \ann_{R}M$, one has 
$\Ext^*_{R_{\fp}}(k(\fp),M_{\fp})\ne 0$. Observe $M_{\fp}\ne 0$ since $\fp$ contains $\ann_{R}M$.
Therefore, localizing at $\fp$ one may assume $R$ is a graded local ring and $M$ is a non-zero finitely generated $R$-module, and then the desired result is that $\Ext^{*}_{R}(k,M)\ne 0$, where $k$ is the graded residue field of $R$. We note that the standard results on associated primes carry over, with identical proofs, to this graded context; one has also an analogue of Nakayama's Lemma. Thus, arguing as in the proof of \cite[(1.2.5)]{BH}, one can deduce the desired non-vanishing.
\end{proof} 

\subsection*{Specialization closed sets}  
Given a subset $\mcU\subseteq\spec R$, we consider the full subcategory of $\Mod R$ with objects
\[
\sfM_\mcU=\{M\in \Mod R\mid \supp_R M\subseteq \mcU\}\,.\]
  
The next result do not hold for arbitrary subsets of $\spec R$. In fact, each of
statements (1) and (2) characterize the property that $\mcV$ is specialization closed.

\begin{lemma}\label{le:specialization}
Let $\mcV$ be a specialization closed subset of $\spec R$.
\begin{enumerate}[{\quad\rm(1)}]
\item
For each $R$-module $M$, one has
\[
\supp_R M\subseteq \mcV\quad\Longleftrightarrow \quad M_\fq=0\text{ for each $\fq$ in
  $\spec R\setminus \mcV$}\,.\]
  \item The subcategory $\sfM_\mcV$ of $\Mod R$ is closed under direct sums, and in any
  exact sequence $0\to M'\to M\to M''\to 0$ of $R$-modules, $M$ is in $\sfM_\mcV$ if and
  only if $M'$ and $M''$ are in $\sfM_\mcV$.\end{enumerate}\end{lemma}

\begin{proof} 
  Since $\mcV$ is specialization closed, it contains $\supp_RM$ if and only if it contains
  $\cl(\supp_RM)$. Thus, (1) is a corollary of Lemma~\ref{le:supp-ann}(1).  Given this (2)
  follows, since for each $\fq$ in $\spec R$, the functor taking an $R$-module $M$ to
  $M_\fq$ is exact and preserves direct sums.
\end{proof}

Torsion modules and local modules can be recognized from their supports:

\begin{lemma} 
\label{le:torsion-local}
Let $\fp$ be a prime ideal in $R$ and let $M$ be an $R$-module. The following statements
hold:
\begin{enumerate}[{\quad\rm(1)}] 
\item $M$ is $\fp$-local if and only if $\supp_R M\subseteq \{\fq\in \spec R\mid
  \fq\subseteq \fp\}$.
\item $M$ is $\fp$-torsion if and only if  $\supp_R M\subseteq\mcV(\fp)$.
\end{enumerate}
\end{lemma}

\begin{proof}
Let $I$ be a minimal injective resolution of $M$.

(1) Since $I_\fp$ is a minimal injective resolution of $M_\fp$, and minimal injective
resolutions are unique up to isomorphism, $M\cong M_\fp$ if and only if $I\cong
I_\fp$. This implies the desired equivalence, by Lemma~\ref{le:inj-local}.

(2) When $\supp_R M\subseteq\mcV(\fp)$, then, by definition of support, one has that $I^0$
is isomorphic to a direct sum of copies of $E(R/\fq)$ with $\fq\in\mcV(\fp)$.  Since each
$E(R/\fq)$ is $\fp$-torsion, so is $I^0$, and hence the same is true of $M$, for it is a
submodule of $I^0$.

Conversely, when $M$ is $\fp$-torsion, $M_\fq =0$ for each $\fq$ in $\spec R$ with
$\fq\not\supseteq\fp$. This implies $\supp_R M\subseteq \mcV(\fp)$, by
Lemma~\ref{le:specialization}(1).
\end{proof} \begin{lemma}\label{le:anti-specialization}
Let $\fp$ be a prime ideal in $R$ and set $\mcU=\{\fq\in\spec R\mid \fq\subseteq \fp\}$.
\begin{enumerate}[{\quad\rm(1)}]
\item The subcategory $\sfM_\mcU$ of $\Mod R$ is closed under taking kernels, cokernels,
  extensions, direct sums, and products.
\item Let $M$ and $N$ be $R$-modules. If $N$ is in $\sfM_\mcU$, then $\Hom^*_R(M,N)$ is in
  $\sfM_\mcU$.\end{enumerate}\end{lemma}

\begin{proof}
  (1) The objects in the subcategory $\sfM_\mcU$ are precisely the
  $\fp$-local $R$-modules, by Lemma~\ref{le:torsion-local}(1).  Thus
  the inclusion functor has a left and a right adjoint. It follows
  that $\sfM_\mcU$ is an abelian full subcategory of $\Mod R$, closed
  under direct sums and products.

  (2) Pick a presentation $F_1\to F_0\to M\to 0$, where $F_0$ and $F_1$ are free
  $R$-modules.  This induces an exact sequence
\[
0\lto \Hom^*_R(M,N)\lto \Hom^*_R(F_0,N)\lto \Hom^*_R(F_1,N)
\]  
Since the $R$-modules $\Hom^*_R(F_i,N)$ are products of shifts of copies of $N$, it
follows from (1) that $\Hom^*_R(M,N)$ is in $\sfM_\mcU$, as claimed.
\end{proof} 

\section{Localization for triangulated categories}\label{se:localization}
Our notion of local cohomology and support for objects in triangulated
categories is based on certain localization functors on the category.
In this section we collect their properties, referring the reader to
Neeman~\cite{Ne} for details.  We should like to emphasize that there
is little in this section which would be unfamiliar to experts, for
the results and arguments have precursors in various contexts; confer
\cite{HPS}, and the work of Alonso Tarrío, Jerem\'ias L\'opez, and
Souto Salorio~\cite{JVS1}. However they have not been written down in
the generality required for our work, so detailed proofs are provided,
if only for our own benefit.

Let $\sfT$ be a triangulated category, and let  $\Si$  denote its suspension functor.

\subsection*{Localization functors}
An exact functor $\bloc\col\sfT\to\sfT$ is called \emph{localization
functor} if there exists a morphism $\eta\col\Id_\sfT\to \bloc$ such
that the morphism $\bloc\eta\col \bloc\to \bloc^2$ is invertible and
$\bloc\eta=\eta \bloc$. Recall that a morphism $\mu\colon F\to G$
between functors is invertible if and only if for each object $X$ the
morphism $\mu{X}\colon FX\to GX$ is an isomorphism.  Note that we only
require the existence of $\eta$; the actual morphism is not part of
the definition of $L$ because it is determined by $L$ up to a unique
isomorphism $L\to L$.

The following lemma provides an alternative description of a
localization functor; it seems to be well known in the context of
monads \cite[Chapter~3]{MacL} but we have no explicit reference.

\begin{lemma}\label{le:loc}
Let $\bloc\col\sfT\to\sfT$ be a functor and $\eta\col\Id_\sfT\to \bloc$  a morphism.
The following conditions are equivalent.  
\begin{enumerate}[{\quad\rm(1)}] 
\item The morphism $\bloc\eta\col \bloc\to \bloc^2$ is invertible and
$\bloc\eta=\eta \bloc$.
\item There exists an adjoint pair of functors $F\col\sfT\to\sfS$ and $G\col\sfS\to\sfT$,
  with $F$ the left adjoint and $G$ the right adjoint, such that $G$ is fully faithful,
  $\bloc=G F$, and $\eta\col\Id_\sfT\to G F$ is the adjunction
  morphism.  
\end{enumerate} 
\end{lemma}

\begin{proof}
(1) $\Rightarrow$ (2): Let $\sfS$ denote the full subcategory of $\sfT$ formed by
  objects $X$ such that $\eta{X}$ is invertible. For each $X\in\sfS$, let $\theta{X}\col
  \bloc{X}\to X$ be the inverse of $\eta{X}$. Define $F\col \sfT\to\sfS$ by
  $F{X}=\bloc{X}$ and let $G\col \sfS\to\sfT$ be the inclusion. It is straightforward to
  check that the maps
\begin{align*}
&\Hom_\sfS(F{X},Y)\lto \Hom_\sfT(X,G{Y}),\quad \alpha\mapsto G{\alpha}\comp\eta{X}\\
&\Hom_\sfT(X,G{Y})\lto\Hom_\sfS(F{X},Y),\quad \beta\mapsto \theta{Y}\comp F{\beta}
\end{align*}
are mutually inverse bijections.  Thus, the functors $F$ and $G$ form an adjoint pair.

(2) $\Rightarrow$ (1): Let $\theta\col FG\to\Id_\sfS$ denote the adjunction morphism. It is
then easily checked that the compositions
\[
F\xra{F\eta}FGF\xra{\theta F} F\quad\text{and}\quad G\xra{\eta G}GFG\xra{G\theta}G\]
  are identity morphisms. Now observe that $\theta$ is invertible because $G$ is fully
faithful.  Therefore $\bloc\eta=GF\eta$ is invertible. Moreover, we have
\[
\bloc\eta=GF\eta=(G\theta F)^{-1}=\eta GF=\eta \bloc\,.\]
  This completes the proof.
\end{proof}

\subsection*{Acyclic and local objects}
Let $\bloc\colon\sfT\to\sfT$ be a localization functor.  An object $X$
in $\sfT$ is said to be \emph{$\bloc$-local} if $\eta{X}$ is an
isomorphism; it is \emph{$\bloc$-acyclic} if $\bloc{X}=0$.  We write
$\Im\bloc$ for the full subcategory of $\sfT$ formed by all
$\bloc$-local objects, and $\Ker\bloc$ for the full subcategory formed
by all $\bloc$-acyclic objects. Note that $\Im\bloc$ equals the
\emph{essential image} of $\bloc$, that is, the full subcategory of
$\sfT$ formed by all objects isomorphic to one of the form $\bloc X$
for some $X$ in $\sfT$.  It is easily checked that $\Im\bloc$ and
$\Ker\bloc$ are triangulated subcategories of $\sfT$. 

Let us mention that $\bloc$ induces an equivalence of categories
$\sfT/{\Ker\bloc}\xra{\sim}\Im\bloc$ where $\sfT/{\Ker\bloc}$ denotes
the Verdier quotient of $\sfT$ with respect to $\Ker\bloc$. This can
be deduced from \cite[I.1.3]{GZ}, but we do not need this fact.

\begin{definition}
\label{def:gamma-exists}  
For each $X$ in $\sfT$, complete the map $\eta{X}$ to an exact
triangle
\[
\gam{X}\xra{\theta{X}}X\xra{\eta{X}}\bloc{X}\lto\,.  
\] 

It follows from the next lemma that one obtains a well defined functor
$\gam\col\sfT\to\sfT$.
\end{definition}

\begin{lemma}\label{le:local-acyclic} 
The functor $\gam$ is exact, and the following properties hold.
\begin{enumerate}[{\quad\rm(1)}]
\item $X\in\sfT$ is $\bloc$-acyclic if and only if $\Hom_\sfT(X,-)=0$ on $\bloc$-local objects;
\item $Y\in\sfT$ is $\bloc$-local if and only if $\Hom_\sfT(-,Y)=0$ on $\bloc$-acyclic objects;
\item  $\gam$ is a right adjoint for the inclusion $\Ker\bloc\to \sfT$;
\item $\bloc$ is a left adjoint for the inclusion $\Im\bloc\to\sfT$.
\end{enumerate}\end{lemma}

\begin{proof}
(1) By Lemma~\ref{le:loc}, one has a factorization $\bloc=GF$, where $G$ is a fully faithful right
adjoint of $F$. In particular, when $X$ is $\bloc$-acyclic, $FX=0$ and then for any $\bloc$-local object
$Y$ one has
\[
\Hom_\sfT(X,Y)\cong\Hom_\sfT(X,GFY)\cong\Hom_\sfT(FX,FY)=0\,.
\]
Conversely, if $X$ in $\sfT$ is such that $\Hom_\sfT(X,Y)=0$ for all $\bloc$-local $Y$, then
\[
\Hom_\sfT(FX,FX)\cong\Hom_\sfT(X,GFX)=0\,,
\]
and hence $FX=0$; that is to say, $X$ is $\bloc$-acyclic.

(2) We have seen in (1) that $\Hom_\sfT(X,Y)=0$ if $X$ is $\bloc$-acyclic and $Y$ is $\bloc$-local. 

Suppose that $Y$ is an object with $\Hom_\sfT(-,Y)=0$ on $\bloc$-acyclic objects. Observe
that $\gam Y$ is $\bloc$-acyclic since $\bloc(\eta Y)$ is an isomorphism and $\bloc$ is exact. 
Thus $\theta Y=0$, so that $\bloc Y\cong Y \oplus \Si \gam{Y}$. Since $\Si\gam{Y}$ is $\bloc$-acyclic
and $\bloc Y$ is $L$-local, the isomorphism implies that $\Si\gam{Y}=0$. This is the desired result.

(3) We noted in (2) that $\gam X$ is $\bloc$-acyclic. For each $\bloc$-acyclic object $W$ the map
$\Hom_\sfT(W,\gam{X})\to\Hom_\sfT(W,X)$ induced by $\theta{X}$ is a bijection since $\Hom_\sfT(W,-)=0$ on $\bloc$-local objects. It follows that the exact triangle in Definition \ref{def:gamma-exists} is unique up to unique isomorphism. Moreover, the assignment $X\mapsto \gam X$ defines a functor, right adjoint to the inclusion functor $\Ker\bloc\to\sfT$.  

The functor $\gam\col \sfT\to\Ker\bloc$ is exact because
it is an adjoint of an exact functor; see \cite[Lemma~5.3.6]{Ne}.

(4) This follows from Lemma~\ref{le:loc}.
\end{proof}   

The functor $\gam\colon\sfT\to\sfT$ is a localization functor for the
opposite category $\sfT^\op$. So we think of $\gam$ as $\bloc$ turned
upside down. Our interpretation of $\gam$ as a local cohomology
functor provides another explanation for using the letter $\gam$; see
Section~\ref{se:commutative}.  We will need to use some rules of
composition for localization functors.

\begin{lemma} \label{le:functors-commute} 
Let $\bloc_1$ and $\bloc_2$ be localization functors for $\sfT$. If each $\bloc_1$-acyclic
object is $\bloc_2$-acyclic, then the following statements hold:
\begin{enumerate}[{\quad\rm(1)}]
\item $\gam_{1}\gam_{2} \cong\gam_{1}\cong\gam_{2}\gam_{1}$ and
  $\bloc_1\bloc_2\cong\bloc_2\cong\bloc_2\bloc_1$.
\item $\gam_{1}\bloc_2= 0 =\bloc_2\gam_{1}$.
\item $\gam_{2}\bloc_1 \cong \bloc_1\gam_{2}$.
\item For each $X$ in $\sfT$, there are exact triangles

\begin{align*}
&\gam_{1}{X}\lto\gam_{2}{X}\lto \gam_{2}\bloc_1{X}\lto \\
&\gam_{2}\bloc_1{X} \lto \bloc_1{X}\lto \bloc_2{X}\lto\,.
\end{align*}\end{enumerate}\end{lemma}

\begin{proof}
  The crucial point is that $\gam_i$ is a right adjoint to the inclusion of
  $\bloc_i$-acyclic objects, and $\bloc_i$ is a left adjoint to the inclusion of
  $\bloc_i$-local objects, see Lemma~\ref{le:local-acyclic}.

(1) is an immediate consequence of our hypothesis. 

(2) follows from (1), since
\[
\gam_{1}\bloc_2\cong\gam_{1}\gam_{2}\bloc_2= 0
    =\bloc_2\bloc_1\gam_{1}\cong\bloc_2\gam_{1}\,.\]
  
(3) and (4): For $X$ in $\sfT$ the hypothesis yields morphisms $\alpha\col \gam_{1}{X}\to
\gam_{2}{X}$ and $\beta\col \bloc_1{X}\to \bloc_2{X}$.  They induce the following
commutative diagram
\[
\xymatrixrowsep{2pc}
\xymatrixcolsep{2pc}
\xymatrix{
\gam_{1}{X}\ar@{->}[r] \ar@{->}[d]^{\alpha} & X \ar@{=}[d]\ar@{->}[r]
                   &\bloc_1{X}\ar@{->}[d]^{\beta}\ar@{->}[r]& \\
\gam_{2}{X}\ar@{->}[r] & X \ar@{->}[r] &\bloc_2{X}\ar@{->}[r]&}\]
  Applying $\gam_{2}$ to the top row and $\bloc_1$ to the bottom row, and bearing in mind
the isomorphisms in (1), one obtains exact triangles
\begin{gather*}
\gam_{1}{X} \lto \gam_{2}{X} \lto \gam_{2}\bloc_1{X}\lto \\
\bloc_1\gam_{2}{X}\lto \bloc_1{X} \lto \bloc_2{X}\lto\,.
\end{gather*}
These exact triangles yield the first and the last isomorphism below:
\[
\gam_{2}\bloc_1{X} \cong \cone(\alpha)\cong \cone(\Si^{-1}\beta)\cong  \bloc_1\gam_{2}{X}\,,\]
  The one in the middle is given by the octahedral axiom.

This completes the proof.
\end{proof}

The example below shows that localization functors usually do not commute.  We are
grateful to Bernhard Keller for suggesting it.

\begin{example}\label{ex:keller}
Let $A=\left[\begin{smallmatrix}k&k\\ 0&k\end{smallmatrix}\right]$ be the algebra of
$2\times 2$ upper triangular matrices over a field $k$ and let $\sfT$ denote the derived
category of all $A$-modules.  Up to isomorphism, there are precisely two indecomposable
projective $A$-modules:
\[
P_1=
\begin{bmatrix} k&0\\    0&0\end{bmatrix} \qquad\text{and}\qquad
P_2= \begin{bmatrix}0&k\\    0&k\end{bmatrix}\]
  satisfying $\Hom_A(P_1,P_2)\neq 0$ and $\Hom_A(P_2,P_1)= 0$.  For $i=1,2$ we let $\bloc_i$
denote the localization functor such that the $\bloc_i$-acyclic objects form the smallest
localizing subcategory containing $P_i$, viewed as a complex concentrated in degree
zero. One then has $\bloc_1\bloc_2\neq\bloc_2\bloc_1$, since
\[
\gam_1\gam_2(P_2)=\gam_1(P_2)=P_1\quad\text{and}\quad
\gam_2\gam_1(P_2)=\gam_2(P_1)=0.\,\]
  \end{example}

\subsection*{Existence} The following criterion for the existence of a localization
functor will be used; it is contained in \cite[Section~7]{Ma}.

\begin{proposition}\label{pr:homloc}
Let $\sfT$ be a triangulated category which admits small coproducts and is compactly
generated, and let $\sfA$ be an abelian Grothendieck category.  Let $H\col\sfT\to\sfA$ be
a cohomological functor which preserves all coproducts.

There then exists a localization functor $\bloc \col\sfT\to\sfT$ with the following
property: For each $X\in\sfT$, one has $\bloc{X}=0$ if and only if $\hh{\Si^nX}=0$ for all
$n\in\bbz$. \qed
\end{proposition} 

\section{Local cohomology}
\label{se:cohomology}
Let $\sfT$ be a compactly generated triangulated category. By this we mean that $\sfT$ is a triangulated category that admits small coproducts, the isomorphism classes of compact objects in $\sfT$ form a set,
and for each $X\in \sfT$, there exists a compact object $C$ such that $\Hom_{T}(C,X)\ne 0$.
Let $\sfT^c$ denote the full subcategory which is formed by all compact objects. We shall identify $\sfT^c$ with a set of representatives for the isomorphism classes of compact objects in $\sfT$ whenever this is convenient.

Let $X,Y$ be objects in $\sfT$. We write $\Hom_{\sfT}(X,Y)$ for the abelian group of morphisms in $\sfT$ from $X$ to $Y$. We consider also the graded abelian group:
\[
\Hom^*_\sfT(X,Y)=\coprod_{i\in\bbz}\Hom_\sfT(X,\Si^i Y)\,.\]
  Set $\End^*_\sfT(X)=\Hom^*_\sfT(X,X)$; it has a natural structure of a graded ring.  The graded abelian group $\Hom^*_\sfT(X,Y)$ is a right-$\End^*_\sfT(X)$ and left-$\End^*_\sfT(Y)$ bimodule.

\subsection*{Center} 
Let $\cent(\sfT)$ denote the graded center of $\sfT$. This is a graded-commutative ring, where, for each $n\in\bbz$, the component in degree $n$ is
\[
\cent(\sfT)^n=\{\eta\col\Id_\sfT\to\Si^n\mid\eta\Si=(-1)^n\Si\eta\}\,.
\]
While $\cent(\sfT)$ may not be a set, this is not an issue for our focus will be on a graded-commutative ring $R$ equipped with a homomorphism $\phi\colon R\to \cent(\sfT)$. What this amounts to is that for each object $X$ in $\sfT$ one has a homomorphism of graded rings
\[
\phi_X\col R\lto\End^*_\sfT(X)\,,\]
  such that the induced actions of $R$ on $\Hom^*_\sfT(X,Y)$, from the right via $\phi_X$ and from the left via $\phi_Y$, are compatible, in the sense that, for any homogeneous elements $r\in R$ and $\alpha\in\Hom^*_\sfT(X,Y)$, one has
\[
\phi_Y(r)\alpha=(-1)^{|r||\alpha|}\alpha\phi_X(r)\,.\]
  In this way, each graded abelian group $\Hom^*_\sfT(X,Y)$ is endowed with a structure of a graded $R$-module. For each $n\in\mathbb Z$ one has a natural isomorphism of $R$-modules:
\[
\Hom^*_\sfT(\Si^n X,Y)=\Hom^*_\sfT(X,Y)[n]\,.\]
  
For example, one has a homomorphism $\mathbb Z\to Z(\sfT)$ sending $n$ to $n\cdot \id \colon \Id_\sfT\to\Si^0$.

\begin{notation}
\label{eqn:daring} 
For the rest of this paper, we fix a graded-commutative noetherian ring $R$ and a homomorphism of graded rings
$R\lto \cent(\sfT)$, and say that $\sfT$ is an \emph{$R$-linear triangulated category}.

Let $C$ be an object in $\sfT$. For each object $X$ in $\sfT$, we set
\[
\coh CX = \Hom^*_{\sfT}(C,X)\,,\]
  and think of this $R$-module as the cohomology of $X$ with respect to $C$.
\end{notation}

The following lemma explains to what extent the cohomology $\coh CX$ of $X$ depends on the choice of the object $C$. A full subcategory of $\sfT$ is \emph{thick} if it is a triangulated subcategory, closed under taking direct summands. Given a set $\sfC$ of objects in $\sfT$, the intersection of all thick subcategories of $\sfT$ containing $\sfC$ is again a thick subcategory, which is said to be \emph{generated} by $\sfC$.

\begin{lemma}\label{le:choice_cohomology}
Let $\sfC$ be a set of objects in $\sfT$ and $C_0$ an object contained in the thick subcategory generated by $\sfC$. Then for each $X\in \sfT$ one has
\[
\supp_R\coh {C_0}X\subseteq\bigcup_{C\in\sfC}\cl(\supp_R\coh CX)\,.\]
  \end{lemma}

\begin{proof}
It suffices to prove that the subcategory of $\sfT$ with objects $D$ such that $\supp_R\coh {D}X$ is a subset of the right hand side of the desired inclusion is thick.  It is clear that it is closed under suspensions and direct summands. Each exact triangle $C'\to C\to C''\to$ induces an exact sequence $\coh {C''}X\to\coh {C}X\to\coh {C'}X$ of $R$-modules, so Lemma~\ref{le:specialization} implies
\[
\supp_R\coh CX\subseteq\cl(\supp_R\coh {C'}X) \cup\cl(\supp_R\coh {C''}X)\,,\]
  which implies that the subcategory is also closed under exact triangles.
\end{proof}

\subsection*{Local cohomology}
In this paragraph we introduce local cohomology with support in specialization closed subsets of $\spec R$. Given a subset $\mcU\subseteq \spec R$, set
\[
\sfT_{\mcU}=\{X\in\sfT\mid \supp_R \coh CX\subseteq \mcU\text{ for each }C\in\sfT^c\}\,.\]
  A subcategory of $\sfT$ is \emph{localizing} if it is thick and closed under taking small coproducts.  Using Lemma~\ref{le:specialization}(2), a routine argument yields the following statement.

\begin{lemma}\label{le:localizing}
If $\mcV\subseteq\spec R$ is specialization closed, then  the  subcategory $\sfT_{\mcV}$ of $\sfT$ is localizing.  \qed
\end{lemma}

A \emph{colocalizing} subcategory is one which is thick and closed under taking small products.  The next result is a direct consequence of Lemma~\ref{le:anti-specialization}(1).

\begin{lemma}\label{le:colocalizing}
Let $\fp$ be a prime ideal in $R$ and set $\mcU(\fp) = \{\fq\in \spec R\mid \fq\subseteq\fp\}$.  The subcategory $\sfT_{\mcU(\fp)}$ of $\sfT$ is localizing and colocalizing.  \qed
\end{lemma}

\begin{proposition}
Let $\mcV\subseteq \spec R$ be specialization closed.  There exists a localization functor $\bloc_{\mcV}\col\sfT\to\sfT$ with the property that $\bloc_{\mcV}X=0$ if and only if $X\in\sfT_{\mcV}$.
\end{proposition}

\begin{proof}
The following functor is cohomological and preserves small coproducts:
\[
\hH \col\sfT\lto\prod_{C\in\sfT^c}\Mod R, \quad \text{where}\quad X\mapsto\big(\coprod_{\fp\not\in\mcV}\coh CX_\fp\big)_{C\in\sfT^c}\,.\]
  Proposition~\ref{pr:homloc} applies and gives a localization functor $\bloc_{\mcV}$ on $\sfT$ with the property that an object $X$ is $\bloc_{\mcV}$ acyclic if and only if $\hh X=0$. It remains to note that this last condition is equivalent to the condition that $X$ is in $\sfT_{\mcV}$, by Lemma~\ref{le:specialization}(1).
\end{proof}

\begin{definition}\label{def:loc-seq}  
Let $\mcV$ be a specialization closed subset of $\spec R$, and $\bloc_{\mcV}$ the associated localization functor given by the proposition above.  By \ref{def:gamma-exists}, one then gets an exact functor $\gam_{\mcV}$ on $\sfT$ and for each object $X$ a natural exact triangle
\[
\gam_{\mcV}X\lto X\lto \bloc_{\mcV}X\lto\,.\]
  We call $\gam_{\mcV}X$ the \emph{local cohomology} of $X$ supported on $\mcV$, and the triangle above the \emph{localization triangle} with respect to $\mcV$.  The terminology may seem unfounded in the first encounter, but look ahead to Theorems~\ref{thm:localization-support} and \ref{thm:ca-main}
\end{definition}

\subsection*{Localization at a prime}
The next result realizes localization on cohomology as a localization on $\sfT$, and explains the nomenclature `localization'. We fix a point $\fp$ in $\spec R$ and set
\[
\mcZ(\fp)=\{\fq\in \spec R\mid \fq\not\subseteq\fp\} \quad\text{and}\quad\mcU(\fp)
=\{\fq\in\spec R\mid \fq\subseteq \fp\}.\]
  
\begin{theorem}\label{thm:locprime}
Let $\fp\in\spec R$. For each compact object $C$ in $\sfT$, the
natural map $X\to \bloc_{\mcZ(\fp)}X$ induces an isomorphism of graded
$R$-modules
\[
\coh CX_\fp\xra{\sim}\coh C{\bloc_{\mcZ(\fp)}X}\,.\]
  \end{theorem}

The proof of this theorem uses an auxiliary construction.

\begin{definition}\label{injectives:lifting}  
Let $I$ be an injective $R$-module and $C$ a compact object in
$\sfT$. Brown representability yields an object $ I_{C}$ in $\sfT$
such that
\[
  \Hom_R(\coh C-,I)\cong\Hom_\sfT(-, I_C)\,.\]
  This isomorphism extends to an isomorphism of functors of graded $R$-modules
\begin{equation}\label{eq:I}
\Hom_R^*(\coh C-,I)\cong\Hom_\sfT^*(-, I_C)\,.
\end{equation}
  \end{definition}

\begin{lemma}\label{le:locp}
Let $\fp\in \spec R$. If $I$ is an injective $R$-module with $\supp_R I\subseteq \mcU(\fp)$, then $I_C\in\sfT_{\mcU(\fp)}$ for each compact object $C$.
\end{lemma}

\begin{proof}
The $R$-modules $\Hom_{\sfT}^*(X, I_C)$ and $\Hom_{R}^*(\coh CX,I)$ are isomorphic for each
$X$ in $\sfT$.  If $\supp_R I\subseteq\mcU(\fp)$, then $\supp_R\Hom_{\sfT}^*(X, I_C)\subseteq \mcU(\fp)$, by Lemma~\ref{le:anti-specialization}(2), so specializing $X$ to compact objects in $\sfT$ yields the desired result.
\end{proof}

\begin{proof}[Proof of Theorem~\emph{\ref{thm:locprime}}]
We note that a compactly generated triangulated category admits small products; see \cite[(8.4.6)]{Ne}.

Let $I$ be the injective hull of $R/\fp$ and $ I_C$ the object in $\sfT$ corresponding to $C$, see
Definition~\ref{injectives:lifting}.  Set $\mcZ=\mcZ(\fp)$ and $\mcU=\mcU(\fp)$.

We claim that $\sfT_{\mcU}$ is perfectly cogenerated in the sense of
\cite[Definition~1]{Kr:br} by the set of objects $I_C$ where $C$ is 
compact in $\sfT$. Indeed, each $ I_C$ belongs to
$\sfT_{\mcU}$, by Lemma~\ref{le:locp}. If $X$ is a nonzero object in
$\sfT_{\mcU}$, then there exists a compact object $C$ such that $\coh CX$ is nonzero;
it is $\fp$-local, by Lemma~\ref{le:torsion-local}, so $\Hom_\sfT^*(X,I_C)\neq 0$, by
\eqref{eq:I}.  Now let $\phi_i\col X_i\to Y_i$ be a family of maps in $\sfT_{\mcU}$ such that
$\Hom_\sfT(Y_i, I_C) \to \Hom_\sfT(X_i, I_C)$ is surjective for all $
I_C$. Then $\coh C{\phi_i}$ is a monomorphism for all $C$ and $i$,
since the injectives $E(R/\fp)[n]$ cogenerate the category of
$\fp$-local $R$-modules. Thus the product $\prod_i\phi_i\col
\prod_iX_i\to \prod_iY_i$ induces a monomorphism
\[
\coh C{\, \prod_i\phi_i\, }= \prod_i\coh C{\phi_i}
\]
and therefore $\Hom_\sfT(\prod_i\phi_i, I_C)$ is surjective for each $C$.

The subcategory $\sfT_{\mcU}$ of $\sfT$ is colocalizing, by
Lemma~\ref{le:colocalizing}, and therefore Brown's representability
theorem \cite[Theorem~A]{Kr:br} yields a left adjoint
$F\col\sfT\to\sfT_{\mcU}$ to the inclusion functor $G\col
\sfT_{\mcU}\to\sfT$. Then $L=GF$ is a localization functor, by
Lemma~\ref{le:loc}. An object $X$ is $\bloc$-acyclic iff $\Hom_\sfT^*(X,I_C)=0$ for all
compact $C$, since the $I_{C}$ cogenerate $\sfT_{\mcU}$. By \eqref{eq:I}, this is the
case iff $\coh CX_\fp= 0$ for all compact $C$; equivalently, iff $X$ is in $\sfT_\mcZ$, by
Lemma~\ref{le:specialization}. Thus the subcategory of $L$-acyclic objects equals $\sfT_{\mcZ}$, so $\bloc$ and $\bloc_{\mcZ}$ coincide. 

The adjunction morphism $\eta{X}\col X\to \bloc_{\mcZ}X$
induces a commutative diagram of natural maps
\[
\xymatrix{
\coh CX \ar[rr]^-{\coh C{\eta{X}}}\ar[d]&&\coh C{\bloc_{\mcZ}X} \ar[d]^{\cong}\\
   \coh CX_\fp\ar[rr]^-{\coh C{\eta{X}}_\fp}_{\cong}&&\coh C{\bloc_{\mcZ}X}_\fp\\ }
  \]
The vertical isomorphism is due to Lemma~\ref{le:torsion-local},
because $\bloc_{\mcZ}{X}$ is in $\sfT_{\mcU}$; the horizontal one
holds because the cone of $\eta{X}$ is in $\sfT_{\mcZ}$.  This
completes the proof.
\end{proof}

The import of the following corollary will be clarified in the next section.

\begin{corollary}\label{cor:gammap}
Let $\fp\in \spec R$ and $X$ an object in $\sfT$.  For each $C\in\sfT^c$, the $R$-module $\coh C{\bloc_{\mcZ(\fp)}\gam_{\mcV(\fp)}X}$ is $\fp$-local and $\fp$-torsion; thus $\bloc_{\mcZ(\fp)}\gam_{\mcV(\fp)}X$ is in
$\sfT_{\{\fp\}}$.
\end{corollary}

\begin{proof}
Set $\gam_{\fp}X= \bloc_{\mcZ(\fp)}\gam_{\mcV(\fp)}X$. The $R$-module $\coh C{\gam_{\fp}X}$ is $\fp$-local, since Theorem~\ref{thm:locprime} yields an isomorphism 
\[
\coh C{\gam_{\fp}X}\cong \coh C{\gam_{\mcV(\fp)}{X}}_\fp\,.\]
  As $\gam_{\mcV(\fp)}{X}$ is in $\sfT_{\mcV(\fp)}$, by Lemma~\ref{le:local-acyclic}, the $R$-module $\coh C{\gam_{\mcV(\fp)}{X}}$ is $\fp$-torsion, by Lemma~\ref{le:torsion-local}, so the same holds of its localization $\coh C{\gam_{\fp}X}$.
\end{proof} 

\section{Support}\label{se:support}
Let $\sfT$ be a compactly generated $R$-linear triangulated category. Recall that $R$ is a graded-commutative noetherian ring; see \ref{eqn:daring}.

\subsection*{Support}\label{support}
Fix a prime ideal $\fp$ in $R$. Set $\mcZ(\fp)=\{\fq\in\spec R\mid \fq\not\subseteq \fp\}$; note that $\mcV(\fp)\setminus\mcZ(\fp)=\{\fp\}$.  We define an exact functor $\gam_{\fp}\col\sfT\to\sfT$ by
\[
\gam_{\fp}X=\bloc_{\mcZ(\fp)}\gam_{\mcV(\fp)}X\quad\text{for each $X\in\sfT$}\,.\]
  
The diagram below displays the natural maps involving the localization functors.
\[
\xymatrix{ \gam_{\mcV(\fp)}X\ar[d]\ar[r]&\gam_{\fp}X\ar[d]\\
X\ar[r]&\bloc_{\mcZ(\fp)}X}
\] 
We define the \emph{support} of an object $X$ in $\sfT$ to be the set
\[
\supp_{R} X=\{\fp\in\spec R\mid\gam_{\fp}X\neq 0\}\subseteq \spec R\,.
\] 
A basic property of supports is immediate from the exactness of the functor $\gam_\fp$: 

\begin{proposition} 
\label{pr:triangle}
\pushQED{\qed}%
For each exact triangle $X\to Y\to Z\to \Si X$ in $\sfT$, one has
\[
\supp_{R} Y\subseteq\supp_{R} X\cup\supp_{R} Z\quad\text{and}\quad\supp_{R}\Si X=\supp_{R} X\,. \qedhere
\]
  \end{proposition} 

The next result relates the support of an object to the support of its cohomology; it is one of the principal results in this work. For a module $M$ over a ring $R$, we write $\Min_R M$ for the set of minimal primes in $\supp_R M$. 

\begin{theorem}
\label{thm:detection}
For each object $X$ in $\sfT$, one has an equality:
\[
  \supp_{R} X = \bigcup_{C\in\sfT^c} \Min_R \coh CX\,.\]
  In particular, $\supp_{R} X=\varnothing$ if and only if $X=0$.
\end{theorem} 

The upshot is that the support of an object can be detected by its cohomology. Observe however that this 
requires one to compute cohomology with respect to all compact objects. It is thus natural to seek efficient ways to compute support. The corollary below is one result in this direction.  Recall that a set of compact objects $\sfG$ is said to \emph{generate} $\sfT$ if it is the smallest localizing subcategory containing $\sfG$. 

\begin{corollary}
\label{cor:support-inclusions} 
If $\sfG$ is a set of compact objects which generate $\sfT$, then
\[
\bigcup_{C\in\sfG} \Min_R \coh CX \subseteq \supp_{R}X \subseteq \bigcup_{C\in\sfG}\cl(\supp_R \coh CX)\,.\]
  \end{corollary}  

\begin{proof} 
Both inclusion follow from Theorem~\ref{thm:detection}; this is clear for the one on the left. For the inclusion on the right one uses in addition Lemma~\ref{le:choice_cohomology}. 
\end{proof} 

\begin{remark}
The inclusions in the preceding corollary can be strict; see Example~\ref{ex:supp-incl}. 
When $\sfT$ can be generated by a single compact object, say $C$, one has that $\supp_{R}X\subseteq \cl(\supp_{R}\coh CX)$. In certain contexts, for example, in that of Section~\ref{se:tensor}, we are able to improve this to $\supp_{R}X\subseteq \supp_{R}\coh CX$, and expect that this holds in general. We will address this problem elsewhere.
\end{remark}

For objects with finite cohomology, support has a transparent expression: 

\begin{theorem}\label{thm:finiteness}
Let $X$ be an object in $\sfT$.
\begin{enumerate}[{\quad\rm(1)}]
\item If $C\in \sfT$ is compact and the $R$-module $\coh CX$ is finitely generated, then
\[
\supp_{R}X \supseteq \supp_R\coh CX\,,\]
  and equality holds when $C$ generates $\sfT$.
\item If $X$ is compact and the $R$-module $\End^*_{\sfT}(X)$ is finitely generated, then
\[
\supp_{R}X = \supp_R\End^*_{\sfT}(X)=\mcV(\ann_R\End^*_{\sfT}(X))\,.\]
  \end{enumerate}
\end{theorem} 

Theorems~\ref{thm:detection} and \ref{thm:finiteness} are proved later in this section.  First, we present some applications.  The next result is a precise expression of the idea that for each specialization closed subset $\mcV$ of $\spec R$ the localization triangle
\[
\gam_{\mcV}X\lto X\lto \bloc_{\mcV}X\lto \]
  separates $X$ into a part supported on $\mcV$ and a part supported on its complement.  For this reason, we henceforth refer to $\gam_{\mcV}X$ as the support of $X$ on $\mcV$, and to $\bloc_{\mcV}X$ as the support of $X$ away from $\mcV$. 

\begin{theorem}\label{thm:localization-support}
Let $\mcV$ be a specialization closed subset of $\spec R$. For each $X$ in $\sfT$ the following equalities hold:
\begin{align*}
&\supp_{R}\gam_{\mcV}X = \mcV\, \cap\, \supp_{R}X\\ 
&\supp_{R}\bloc_{\mcV}X = \big(\spec R\setminus \mcV\big)\, \cap\, \supp_{R}X\,.
\end{align*}
\end{theorem} 

\begin{proof}
For each $\fp$ in $\mcV$, with $\mcZ(\fp)=\{\fq\in\spec R\mid \fq\not\subseteq \fp\}$, one has that
\[
\gam_{\fp}\bloc_{\mcV}X=\bloc_{\mcZ(\fp)}\gam_{\mcV(\fp)}\bloc_{\mcV}X =0\,.
\] 
The second equality holds by Lemma~\ref{le:functors-commute}, because $\mcV(\fp)\subseteq \mcV$.  Therefore
\[
\supp_{R}\bloc_{\mcV}X\subseteq (\spec R)\setminus \mcV\,.\]
  On the other hand, since $\gam_{\mcV}X$ is in $\sfT_{\mcV}$, Theorem~\ref{thm:detection} yields an inclusion:
\[
\supp_{R}\gam_{\mcV}X\subseteq \mcV\,.\]
  The desired result readily follows by combining the inclusions above and Proposition~\ref{pr:triangle}, applied to the exact triangle $\gam_{\mcV}X\to X\to\bloc_{\mcV}X\to$.
\end{proof} 

This theorem has a direct corollary, in view of the localization triangle of $\mcV$. 

\begin{corollary}\label{cor:support}  
Let $\mcV$ be a specialization closed subset and $X$ an object in $\sfT$. The following statements hold.
\begin{enumerate}[{\quad\rm(1)}]
\item $\supp_{R}X\subseteq \mcV$ if and only if $X\in\sfT_\mcV$, if and only if the natural map $\gam_{\mcV}X\to X$ is an isomorphism, if and only if $\bloc_{\mcV}X=0$;
\item
$\mcV\cap \supp_{R}X =\varnothing$ if and only if the natural map $X\to \bloc_{\mcV}X$ is an isomorphism, if and only if $\gam_{\mcV}X=0$. \qed\end{enumerate}\end{corollary} 

\begin{corollary}\label{cor:orthogonality}
If $X$ and $Y$ are objects in $\sfT$ such that $\cl(\supp_{R} X)\cap\supp_{R} Y=\varnothing$, then $\Hom_\sfT(X,Y)=0$.
\end{corollary} 

\begin{proof}
Set $\mcV=\cl(\supp_{R}X)$. Then $X$ is $\bloc_\mcV$-acyclic, while $Y$ is $\bloc_\mcV$-local, by Corollary~\ref{cor:support}. It follows from Lemma~\ref{le:local-acyclic} that $\Hom_\sfT(X,Y)=0$.
\end{proof} 

The next result builds also on Corollary~\ref{cor:gammap}. 

\begin{corollary}
Let $\fp$ be a point in $\spec R$ and $X$ a nonzero object in $\sfT$. The following conditions are equivalent.
\begin{enumerate}[\quad\rm(1)]
\item  $\gam_{\fp}X\cong X$;
\item  $\supp_RX=\{\fp\}$;
\item  $X\in \sfT_{\{\fp\}}$.\end{enumerate}\end{corollary} 

\begin{proof}
Corollary~\ref{cor:gammap} yields that (1) implies (3), while Theorem~\ref{thm:detection} yields that (3) implies (2). 

(2) $\implies$ (1): As $\supp_RX=\{\fp\}$ Theorem~\ref{thm:localization-support} implies $\supp_R\gam_{\mcV(\fp)}X=\{\fp\}$ as well. Therefore, Corollary~\ref{cor:support} yields that the natural maps
\[
\gam_{\mcV(\fp)}X\to X\quad\text{and}\quad 
\gam_{\mcV(\fp)}X\to \bloc_{\mcZ(\fp)}\gam_{\mcV(\fp)}X =\gam_{\fp}X\]
  are both isomorphisms. Thus, (1) holds.
\end{proof} 

The proofs of Theorems~\ref{thm:detection} and \ref{thm:finiteness} involve `Koszul objects' discussed in the following paragraphs, which include auxiliary results of independent interest.  

\begin{definition}
Let $r$ be an element in $R$; it is assumed to be homogeneous, since we are in the category of graded $R$-modules. Set $d=|r|$, the degree of $r$. Let $C$ be an object in $\sfT$.  We denote by $\kos Cr$ any object that appears in an exact triangle
\[
C\stackrel{r}\lto \Si^{d}C\lto \kos Cr \lto \,,\]
  and call it a \emph{Koszul object of $r$ on $C$}; it is well defined up to (nonunique) isomorphism. For any object $X$ in $\sfT$, applying $\Hom^*_{\sfT}(-,X)$ to the triangle above yields an exact sequence of
$R$-modules:
\begin{multline}\label{eq:koszul-les}
\Hom^*_{\sfT}(C,X)[d+1] \stackrel{\mp r}\lto\Hom^*_{\sfT}(C,X)[1]\lto\\
\lto \Hom^*_{\sfT}(\kos Cr,X)\lto \Hom^*_{\sfT}(C,X)[d]\stackrel{\pm r}\lto\Hom^*_{\sfT}(C,X)
\end{multline}
Applying the functor $\Hom^*_{\sfT}(X,-)$ results in a similar exact sequence. Given a sequence of elements $\bfr=r_1,\ldots,r_n$ in $R$, consider objects $C_i$ defined by
\begin{equation}\label{eq:koszul}  
C_i = \begin{cases}
C & \text{for $i=0$,}\\
\kos{C_{i-1}}{r_i} & \text{for $i\geq 1$.}
\end{cases}
\end{equation}
Set $\kos C{\bfr} = C_n$; this is a \emph{Koszul object of $\bfr$ on $C$}.  Finally, given an ideal $\fa$ in $R$, we write $\kos C{\fa}$ for any Koszul object on $C$, with respect to some finite sequence of generators for $\fa$. This object may depend on the choice of the minimal generating sequence for $\fa$. Note that when $C$ is compact, so is $\kos C{\fa}$.
\end{definition} 

\begin{lemma}\label{le:koszul-properties}
Let $\fp$ be a point in $\spec R$, let $C$ be an object in $\sfT$, and $\kos C{\fp}$ a Koszul object of $\fp$ on $C$.  For each object $X$ in $\sfT$, the following statements hold.
\begin{enumerate}[{\quad\rm(1)}]
\item There exists an integer $s\ge 0$ such that 
\[
\fp^s \Hom_\sfT^*(X,\kos C{\fp})=0=\fp^s \Hom_\sfT^*(\kos C{\fp},X)\,.\]
  Therefore, the  $R$-modules $\Hom_\sfT^*(X,\kos C{\fp})$ and $\Hom_\sfT^*(\kos C{\fp},X)$ are $\fp$-torsion.
\item
The Koszul object $\kos C{\fp}$ is in $\sfT_{\mcV(\fp)}$.
\item 
If $\Hom^*_{\sfT}(C,X)=0$, then $\Hom_\sfT^*(\kos C{\fp},X)=0$. The converse holds when the $R$-module $\Hom^*_{\sfT}(C,X)$ is $\fp$-torsion, or $\fp$-local and finitely generated over $R_\fp$.
\item When $C$ is compact, one has
\[
\Hom^*_\sfT(\kos C{\fp},\gam_{\fp}X) \cong \Hom^*_\sfT(\kos C{\fp},X)_\fp\,.\]
  \end{enumerate}
\end{lemma} 

\begin{proof}
By definition, $\kos C{\fp}$ is obtained as an iterated cone on a sequence $r_1,\dots,r_n$, generating the ideal $\fp$.  Let $C_i$ be the objects defined in \eqref{eq:koszul}; thus, $C_n=\kos C{\fp}$.
 
  (1) We construct such an integer $s$ by an iteration. Set $s_0=1$ and assume that for some $i\ge 0$ there exists an $s_i$ with 
\[
  (r_1,\dots,r_i)^{s_i}\Hom_\sfT^*(C_i,X)=0\,.\]
  Since $C_{i+1}$ is the cone of the morphism $C_i\xra{r_{i+1}}C_i$, it follows from \eqref{eq:koszul-les} that $(r_{i+1}^2)$ and $(r_1,\dots,r_i)^{2s_i}$ both annihilate $\Hom_\sfT^*(C_{i+1},X)$.  Thus, set $s_{i+1}=2s_i+1$, and repeat the process.
  
A similar argument proves the claim about $\Hom_\sfT^*(X,\kos C{\fp})$. 

(2) follows from (1), by Lemma~\ref{le:torsion-local}. 

(3) If $\Hom^*_{\sfT}(C,X)=0$, repeated application of \eqref{eq:koszul-les} yields $\Hom^*_{\sfT}(C_n,X)=0$. 
Suppose $\Hom^*_{\sfT}(C,X)\ne 0$. When the $R$-module $\Hom^*_{\sfT}(C,X)$ is $\fp$-torsion, it follows from from the exact sequence \eqref{eq:koszul-les} that $\Hom^*_{\sfT}(C_1,X)$ is again $\fp$-torsion; that it
is also nonzero is immediate from the same sequence for $r_1$ is in $\fp$.  An iteration yields $\Hom^*_{\sfT}(C_n,X)\ne 0$. 

When $\Hom^*_{\sfT}(C,X)$ is $\fp$-local and finitely generated over $R_\fp$, the exact sequence \eqref{eq:koszul-les} implies the same is true of $\Hom^*_{\sfT}(C_1,X)$. Nakayama's lemma implies it is nonzero since $r_1$ is in $\fp$. Again, an iteration yields the desired nonvanishing. 

(4) Since $\kos C{\fp}$ is in $\sfT_{\mcV(\fp)}$, by (2), the map $\gam_{\mcV(\fp)}{X}\to X$ induces an isomorphism
\[
\Hom^*_\sfT(\kos C{\fp},\gam_{\mcV(\fp)}{X})\cong \Hom^*_\sfT(\kos C{\fp},X)\,.\]
  With $\mcZ(\fp)=\{\fq\in\spec R\mid \fq\not\subseteq\fq\}$, the morphism
\[
\gam_{\mcV(\fp)}{X}\lto \bloc_{\mcZ(\fp)}\gam_{\mcV(\fp)}{X}=\gam_{\fp}X\]
  and Theorem~\ref{thm:locprime} induce an isomorphism
\[
\Hom^*_\sfT(\kos C{\fp},\gam_{\mcV(\fp)}{X})_\fp\cong
\Hom^*_\sfT(\kos C{\fp},\gam_{\fp}X)\,.
\] 
The desired isomorphism follows by combining both isomorphisms.
\end{proof} 

\subsection*{Koszul objects and support}
Next we express the support of an object in $\sfT$ via Koszul objects; this leads to a proof of Theorem~\ref{thm:detection}.  

\begin{proposition}\label{pr:koszul-support}
Let $X$ be an object in $\sfT$.  For each point $\fp$ in $\spec R$ and compact object $C$ in $\sfT$, the following conditions are equivalent:
\begin{enumerate}[{\quad\rm(1)}]
\item $\Hom_\sfT^*(C,\gam_{\fp}X)\neq 0$;
\item $\Hom_\sfT^*(\kos C{\fp},\gam_{\fp}X)\neq 0$;
\item $\Hom^*_\sfT(\kos C{\fp},X)_\fp\neq 0$;
\item $\fp\in\supp_R\Hom^*_\sfT(\kos C{\fp},X)$.\end{enumerate}\end{proposition}

\begin{proof}
  (1) $\iff$ (2): This is a consequence of Lemma~\ref{le:koszul-properties}(3), for $\Hom_\sfT^*(C,\gam_\fp X)$ is $\fp$-torsion, by Corollary~\ref{cor:gammap}. 

(2) $\iff$ (3) follows from Lemma~\ref{le:koszul-properties}(4). 

(3) $\iff$ (4): Set $M=\Hom^*_\sfT(\kos C{\fp},X)$. By Lemma~\ref{le:koszul-properties}(1), the $R$-module $M$ is $\fp$-torsion, so $\supp_RM\subseteq \mcV(\fp)$. Therefore, Lemma~\ref{le:supp-ann}(2) yields:
\[
\supp_{R}(M_\fp)= \supp_R{M}\cap\{\fp\}\,.\]
  This implies the desired equivalence.
\end{proof} 

We are now in a position to state a prove a more refined version of Theorem~\ref{thm:detection}.
Compare it with Corollary~\ref{cor:support-inclusions}.

\begin{theorem}
\label{thm:detectionplus}
Let $\sfG$ be a set of compact generators for $\sfT$. For each object $X$ in $\sfT$, one has the following equalities:
\[
  \bigcup_{\substack{C\in \sfG \\\fp\in\spec R}} \Min_R\coh {\kos C{\fp}}X
  = \supp_{R} X = \bigcup_{C\in\sfT^c} \Min_R \coh CX \,.
\]
In particular, $\supp_{R} X=\varnothing$ if and only if $X=0$.
\end{theorem}

We prove this result concurrently with Theorem~\ref{thm:finiteness}.

\begin{proof}[Proofs of Theorems~\emph{\ref{thm:detectionplus}} and \emph{\ref{thm:finiteness}}]
First we verify that each object $X$ in $\sfT$ satisfies:
\[
\supp_{R} X \subseteq \bigcup_{\substack{C\in \sfG \\\fp\in\spec R}} \Min_R\coh {\kos C{\fp}}X \subseteq
\bigcup_{C\in\sfT^c} \Min_R \coh CX\,.
\]
When $C$ is compact, so is $\kos C{\fp}$, hence the inclusion on the right is obvious. 

Let $\fp$ be a point $\spec R$ such that $\gam_{\fp}X\ne 0$. Since $\sfG$ generate $\sfT$, there exists an object $C$ in $\sfG$ with $\Hom^*_{\sfT}(C,{\gam_{\fp}X})\ne 0$. Proposition~\ref{pr:koszul-support} yields that $\fp$ is in the support of the $R$-module $\Hom^*_{\sfT}(\kos C{\fp},X)$. By Lemma~\ref{le:koszul-properties}(1) the latter module is also $\fp$-torsion, so Lemma~\ref{le:torsion-local}(2) implies that its support is contained in $\mcV(\fp)$. Therefore, one obtains that $\fp$ is the minimal prime in the support of the $R$-module $\Hom^*_{\sfT}(\kos C{\fp},X)$.  This justifies the inclusion on the right. 

Let now $C$ be a compact object in $\sfT$ and $\fp$ a point in $\supp_R\Hom^*_{\sfT}(C,X)$. Recall Theorem~\ref{thm:locprime}: for any compact object $D$ in $\sfT$ one has an isomorphism of $R$-modules
\[
\Hom^*_{\sfT}(D,\bloc_{\mcZ(\fp)}{X})\cong \Hom^*_{\sfT}(D,X)_\fp\,.
\]
When $\fp$ is minimal in the support of $\Hom^*_{\sfT}(C,X)$, the $R$-module ${\Hom^*_{\sfT}(C,X)}_\fp$ is nonzero and $\fp$-torsion; when the $R$-module $\Hom^*_{\sfT}(C,X)$ is finitely generated, so is the $R_\fp$-module $\Hom^*_{\sfT}(C,X)_\fp$. Thus, in either case the isomorphism above, applied with $D=C$, and Lemma~\ref{le:koszul-properties}(3) imply that $\Hom^*_{\sfT}(\kos C{\fp},\bloc_{\mcZ(\fp)}{X})$ is nonzero.  Hence the isomorphism above, now applied with $D=\kos C{\fp}$, and Lemma~\ref{le:koszul-properties}(4) yields that $\gam_{\fp}X$ is nonzero, that is to say, $\fp$ is in $\supp_{R} X$. 

At this point, we have proved Theorem~\ref{thm:detectionplus}, and the first claim in Theorem~\ref{thm:finiteness}. In case that $C$ generates $\sfT$, the inclusion $\supp_{R} X\subseteq\supp_R\coh CX$ follows from Corollary~\ref{cor:support-inclusions} because $\supp_R\coh CX$ is specialization closed, by  Lemma~\ref{le:supp-ann}(1). 

In the remainder of the proof, $X$ is assumed to be compact and the $R$-module
$\End^*_{\sfT}(X)$ is finitely generated.  Set $\fa=\ann_R\End^*_{\sfT}(X)$. For each $C$ in $\sfT^c$, the $R$-action on $\coh CX$ factors through the homomorphism of rings $R\to R/\fa$. Therefore Lemma~\ref{le:supp-ann}(1) yields
\[
\supp_R\coh CX \subseteq \mcV(\fa)\,.\]
  In view of Theorem~\ref{thm:detection}, which has been proved, this yields the first inclusion below:
\[
\supp_{R}X \subseteq \mcV(\fa) = \supp_R\End^*_\sfT(X)\subseteq \supp_{R}X\,.\]
  The equality holds by Lemma~\ref{le:supp-ann}(1) and the last inclusion holds by the already established part of Theorem~\ref{thm:finiteness}; to invoke either result, one requires the hypothesis that the $R$-module $\End^*_\sfT(X)$ is finitely generated.  This completes the proof of Theorem~\ref{thm:finiteness} as well.
\end{proof} 

\begin{remark}
\label{rem:finiteness-refinement}
The proof of Theorem~\ref{thm:finiteness} can be adapted to establish an inclusion
\[
\supp_R\coh CX\subseteq \supp_{R}X
\]  
where $C$ is  compact and the $R$-module $\coh CX$ is finitely generated, under other conditions as well; for instance, when the ring $R$ is concentrated in degree $0$.
\end{remark} 

Supports can be characterized by four, entirely reasonable, properties. 

\subsection*{Axiomatic characterization of support}
Let $\sfG$ be a set of compact generators for $\sfT$. For simplicity, for each $X$ in $\sfT$ set
\[
\supp_R H^*(X)=\bigcup_{C\in\sfG} \supp_R \coh CX\,.\]
  It follows from Lemma~\ref{le:choice_cohomology} that $\cl(\supp_R H^*(X))$ is independent of the choice of a generating set $\sfG$. The result below contains the first theorem in the introduction. 

\begin{theorem}\label{thm:axioms}
There exists a unique assignment sending each object $X$ in $\sfT$ to a subset $\supp_{R} X$ of $\spec R$ such that the following properties hold:
\begin{enumerate}[{\quad\rm(1)}]
\item \emph{Cohomology:} For each object $X$ in $\sfT$ one has
\[
\cl(\supp_{R} X) = \cl(\supp_R H^*(X))\,.
\] 
\item \emph{Orthogonality:} For objects $X$ and $Y$ in $\sfT$, one has that
\[
\cl(\supp_{R} X)\cap\supp_{R} Y=\varnothing \quad\text{implies}\quad \Hom_\sfT(X,Y)=0\,.
\]
  \item \emph{Exactness:} For every exact triangle $W\to X\to Y\to$ in $\sfT$, one has
\[
\supp_{R} X\subseteq \supp_{R} W \cup \supp_{R} Y\,.
\] 
\item \emph{Separation:} For any specialization closed subset $\mcV$ of $\spec R$ and object $X$ in $\sfT$, there exists an exact triangle $X'\to X\to X''\to$ in $\sfT$ such that
\[
\supp_{R} X'\subseteq\mcV \quad\text{and}
\quad\supp_{R}X'' \subseteq\spec R\setminus \mcV\,.
\]
  \end{enumerate}
\end{theorem} 

\begin{proof}
 Corollary~\ref{cor:support-inclusions} implies (1), Corollary~\ref{cor:orthogonality} is (2), Proposition~\ref{pr:triangle} is (3), and, given the localization triangle~\ref{def:loc-seq}, Theorem~\ref{thm:localization-support} entails (4). 
 
Now let $\si\col \sfT\to \spec R$ be a map satisfying properties (1)--(4). 

Fix a specialization closed subset $\mcV\subseteq\spec R$ and an object $X\in \sfT$. 
It suffices to verify that the following equalities hold:
\[
\si(\gam_{\mcV}X)= \si(X)\cap \mcV \quad\text{and}\quad
\si(\bloc_{\mcV}X)= \si(X)\cap(\spec R\setminus \mcV)\,.\tag{$\ast$}\]
  Indeed, for any point $\fp$ in $\spec R$ one then obtains that
\begin{align*}
\si(\gam_{\fp}X) &=\si(\bloc_{\mcZ(\fp)}\gam_{\mcV(\fp)}{X})   \\
                     &=\si(\gam_{\mcV(\fp)}{X})\cap (\spec R\setminus \mcZ(\fp)) \\
                     &=\si(X)\cap \mcV(\fp) \cap (\spec R\setminus \mcZ(\fp))\\
                     &=\si(X)\cap \{\fp\}\,.
\end{align*}
Therefore, $\fp\in \si(X)$ if and only if $\si(\gam_{\fp}X)\ne \varnothing$; this last condition is equivalent to $\gam_{\fp}X\ne 0$, by the cohomology property. The upshot is that $\fp\in \si(X)$ if and only if $\fp\in \supp_{R}X$, which is the desired conclusion. 

It thus remains to prove ($\ast$). 

Let $X'\to X\to X''\to$ be the triangle associated to $\mcV$, provided by property (4). It suffices to verify the following statements:
\begin{enumerate}[{\quad\rm(i)}]
\item 
$\si(X')= \si(X)\cap \mcV$ and $\si(X'')= \si(X)\cap(\spec R\setminus \mcV)$;
\item
$\gam_{\mcV}X\cong X'$ and $\bloc_{\mcV}X\cong X''$.
\end{enumerate} 

The equalities in (i) are immediate from properties (3) and (4). In verifying (ii), the crucial observation is that, by the cohomology property, for any $Y$ in $\sfT$ one has
\[
\si(Y)\subseteq \mcV\iff Y\in\sfT_\mcV\,.\]
  Thus $X'$ is $L_\mcV$-acyclic. On the other hand, property (2) and Lemma~\ref{le:local-acyclic} imply that $X''$ is $\bloc_\mcV$-local.  One thus obtains the following morphism of triangles
\[
\xymatrixrowsep{2pc}
\xymatrixcolsep{2pc}
\xymatrix{
X'\ar@{->}[r] \ar@{->}[d]^{\alpha} & X \ar@{=}[d]\ar@{->}[r]
                   &X''\ar@{->}[d]^{\beta}\ar@{->}[r]& \\
\gam_{\mcV}{X}\ar@{->}[r] & X \ar@{->}[r] &\bloc_{\mcV}{X}\ar@{->}[r]&}
\]
  where the object $\cone(\alpha)\cong\cone(\Si^{-1}\beta)$ is $\bloc_\mcV$-acyclic and $\bloc_\mcV$-local, hence trivial. Therefore, $\alpha$ and $\beta$ are isomorphisms, which yields (ii). 

This completes the proof of the theorem.
\end{proof} 

\begin{remark}
In Theorem~\ref{thm:axioms}, when proving that any function $\sigma$
with properties (1)--(4) coincides with support, properties (3) and
(4) were used only to obtain an exact triangle $X'\to X\to X''\to$
satisfying conditions (i) and (ii), in the proof of the theorem. One may thus replace those
properties by the following one:
\begin{enumerate}[{\quad\rm(1)}]
\item[(3$'$)] Exact separation: \emph{For any specialization closed
subset $\mcV$ of $\spec R$ and object $X$ in $\sfT$, there exists an
exact triangle $X'\to X\to X''\to$ in $\sfT$ such that
\[
\supp_{R} X' = \supp_{R}X\cap\mcV \quad\text{and}
\quad\supp_{R}X'' = \supp_{R}X\cap(\spec R\setminus \mcV)\,.
\]
}
\end{enumerate}
\end{remark}

\section{Properties of local cohomology}\label{se:local cohomology properties} 
Let $\sfT$ be a compactly generated $R$-linear triangulated category; see \ref{eqn:daring}. 

\subsection*{Composition laws}
We provide commutation rules for local cohomology and localization functors and give alternative descriptions of the functor $\gam_{\fp}$. Compare the result below with Example~\ref{ex:keller}. 

\begin{proposition}\label{pr:localizations-commute}
Let $\mcV$ and $\mcW$ be specialization closed subsets of $\spec R$. There are natural isomorphisms of functors
\begin{enumerate}[{\quad\rm(1)}]
\item $\gam_{\mcV}\gam_{\mcW}\cong \gam_{\mcV\cap\mcW} \cong \gam_{\mcW}\gam_{\mcV}$;
\item $\bloc_{\mcV}\bloc_{\mcW}\cong \bloc_{\mcV\cup\mcW} \cong \bloc_{\mcW}\bloc_{\mcV}$;
\item $\gam_{\mcV}\bloc_{\mcW}\cong \bloc_{\mcW}\gam_{\mcV}$.\end{enumerate}\end{proposition} 

\begin{proof}
(1) It suffices to verify the isomorphism on the left.  Let $X$ be an object in $\sfT$, and consider the localization triangle
\[
\gam_{\mcV\cap\mcW}\gam_{\mcV}\gam_{\mcW}X\lto \gam_{\mcV}\gam_{\mcW}X
\lto \bloc_{\mcV\cap\mcW} \gam_{\mcV}\gam_{\mcW}X\lto\]
  of $\gam_{\mcV}\gam_{\mcW}X$ with respect to $\mcV\cap\mcW$. Lemma~\ref{le:functors-commute} provides isomorphisms 
\[
\gam_{\mcV\cap\mcW}\gam_{\mcV}\gam_{\mcW}X \cong
 \gam_{\mcV\cap\mcW}\gam_{\mcW}X\cong  
      \gam_{\mcV\cap\mcW}X\,.\]
  It follows by Theorem~\ref{thm:localization-support} that $\supp_{R}\gam_{\mcV}\gam_{\mcW}X\subseteq \mcV\cap \mcW$, and hence by Corollary~\ref{cor:support} that $\gam_{\mcV}\gam_{\mcW}X$ is in $\sfT_{\mcV\cap\mcW}$.
The exact triangle above now yields the desired result. 

(2) The proof is similar to that of (1). 

(3) The inclusion $\mcV\subseteq \mcV\cup\mcW$ induces a morphism $\theta\col \gam_{\mcV}\to \gam_{\mcV\cup\mcW}$.  We claim that the induced morphisms
\[
\bloc_{\mcW}\theta\col \bloc_{\mcW}\gam_{\mcV}\lto \bloc_{\mcW}\gam_{\mcV\cup\mcW}
\quad\text{and}\quad
\theta\bloc_{\mcW}\col \gam_{\mcV}\bloc_{\mcW} \lto \gam_{\mcV\cup\mcW}\bloc_{\mcW}\]
  are isomorphisms. This implies the desired isomorphism, since Lemma~\ref{le:functors-commute} yields
\[
\bloc_{\mcW}\gam_{\mcV\cup\mcW}\cong\gam_{\mcV\cup\mcW}\bloc_{\mcW}\,.\]
  To prove the claim, consider for each $Y$ in $\sfT$ the exact triangle
\[
\gam_{\mcV}Y \xra{\theta Y} \gam_{\mcV\cup\mcW}Y \lto
\bloc_{\mcV}\gam_{\mcV\cup\mcW}Y\lto\,.\]
  It remains to note that one has isomorphisms
\begin{gather*}
\bloc_{\mcW}\bloc_{\mcV}\gam_{\mcV\cup\mcW}X\cong
\bloc_{\mcW\cup\mcV}\gam_{\mcV\cup\mcW}X =0\\
\bloc_{\mcV}\gam_{\mcV\cup\mcW}\bloc_{\mcW}X\cong
\bloc_{\mcV}\bloc_{\mcW} \gam_{\mcV\cup\mcW}X \cong
\bloc_{\mcW\cup\mcV} \gam_{\mcV\cup\mcW}X=0\,,
\end{gather*}
where the first isomorphism in the first row and the second one in the second row hold by part (2). The first isomorphism in the second row holds by Lemma~\ref{le:functors-commute}.
\end{proof} 

\begin{theorem}\label{thm:gammap-variants}
Let $\fp$ be a point in $\spec R$, and let $\mcV$ and $\mcW$ be specialization closed subsets of $\spec R$ such that $\mcV\setminus\mcW=\{\fp\}$. There are natural isomorphisms
\[
 \bloc_{\mcW}\gam_{\mcV}\cong \gam_\fp\cong \gam_{\mcV}\bloc_{\mcW}\,.\]
  \end{theorem} 

Note that the hypothesis is equivalent to: $\mcV\not\subseteq \mcW$ and $\mcV\subseteq \mcW\cup\{\fp\}$. Our choice of notation is illustrated by the following diagram.
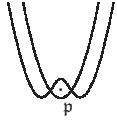
\[
\xy 
(-1,24);(13,24) **\crv{(3,4)&(6,24)&(9,4)}; (1,24);(11,24) **\crv{(6,-1.2)};
(6,12.5)*{.};(7,10)*{\scriptstyle\mathfrak p};
\endxy\]
  \vspace{-27pt} 

\begin{proof}
Set $\mcZ(\fp)=\{\fq\in\spec R\mid \fq\not\subseteq \fp\}$ and $\mcY(\fp)=\mcZ(\fp)\cup\{\fp\}$. The hypothesis implies the following inclusions
\[
\mcV(\fp)\subseteq \mcV\subseteq \mcW\cup \{\fp\}\,,\quad \mcW\subseteq \mcZ(\fp)\,,\quad
\text{and}\quad \mcZ(\fp)\cap \mcV(\fp) \subseteq \mcW\,.\]
  The isomorphisms of localization functors associated to these subsets are used without further comment, see Lemma~\ref{le:functors-commute} and Proposition~\ref{pr:localizations-commute}. 

Given that localization functors commute, it is enough to prove that there are natural isomorphisms
\[
\gam_{\fp} = \bloc_{\mcZ(\fp)}\gam_{\mcV(\fp)}\cong 
\bloc_{\mcW}\gam_{\mcV(\fp)} \cong \bloc_{\mcW}\gam_{\mcV}\,.
\]  

Let $X$ be an object in $\sfT$. We consider the following exact triangles.
\begin{gather*}
\gam_{\mcZ(\fp)}\bloc_{\mcW}\gam_{\mcV(\fp)}X\lto \bloc_{\mcW}\gam_{\mcV(\fp)}X \lto 
\bloc_{\mcZ(\fp)}\gam_{\mcV(\fp)}X\lto\\
\bloc_{\mcW}\gam_{\mcV(\fp)}X\lto \bloc_{\mcW}\gam_{\mcV}X \lto 
\bloc_{\mcW}\bloc_{\mcV(\fp)}\gam_{\mcV}X\lto\,.
\end{gather*}
The first triangle is the localization triangle of the object $\bloc_{\mcW}\gam_{\mcV(\fp)}X$ with respect to $\mcZ(\fp)$. The second triangle is obtained by applying the functor $\bloc_{\mcW}\gam_{\mcV}$ to the localization triangle of $X$ with respect to $\mcV(\fp)$.  It remains to note the isomorphisms:
\begin{gather*}
\gam_{\mcZ(\fp)}\bloc_{\mcW}\gam_{\mcV(\fp)} \cong
      \bloc_{\mcW}\gam_{\mcZ(\fp)}\gam_{\mcV(\fp)} \cong
      \bloc_{\mcW}\gam_{\mcZ(\fp)\cap \mcV(\fp)} = 0\\
      \bloc_{\mcW}\bloc_{\mcV(\fp)}\gam_{\mcV}\cong
      \bloc_{\mcW\cup\mcV(\fp)}\gam_{\mcV}= 0\,.
\end{gather*}
This completes the proof of the theorem.
\end{proof} 

\begin{corollary}\label{cor:gammap-properties}
Let $\fp$ be a point in $\spec R$, and set $\mcZ(\fp)=\{\fq\in\spec R\mid \fq\not\subseteq \fp\}$ and $\mcY(\fp)= \mcZ(\fp)\cup\{\fp\}$.  There are natural isomorphisms
\[
\gam_\fp\cong\gam_{\mcV(\fp)}\bloc_{\mcZ(\fp)}\cong\gam_{\mcY(\fp)}\bloc_{\mcZ(\fp)}\cong
\bloc_{\mcZ(\fp)}\gam_{\mcY(\fp)}\,.\]
  Moreover, for each $X$ in $\sfT$, there are exact triangles
\begin{align*}
&\gam_{\mcZ(\fp)}X\lto\gam_{\mcY(\fp)}X\lto\gam_{\fp}X\lto \\
&\gam_{\fp}X\lto \bloc_{\mcZ(\fp)}X\lto \bloc_{\mcY(\fp)}X\lto\,.
\end{align*}
\end{corollary} 

\begin{proof}
  The isomorphisms are all special cases of Theorem~\ref{thm:gammap-variants}. Given these, the exact triangles follow from Lemma~\ref{le:functors-commute}, since $\mcZ(\fp)\subset \mcY(\fp)$.
\end{proof}  

\subsection*{Smashing localization}\label{Smashing localization}
Let $\sfC$ be a class of objects in $\sfT$; we write $\Loc(\sfC)$ for the smallest localizing subcategory of $\sfT$ containing it.  By definition, $\sfC$ generates $\sfT$ if and only if $\Loc(\sfC)=\sfT$.  

\begin{theorem}\label{thm:localization-finite}
Let $\mcV\subseteq \spec R$ be a specialization closed subset of $\spec R$, and let $\sfG$ be a set of compact generators for $\sfT$. One then has
\[
\sfT_{\mcV}= \Loc(\sfT^c\cap \sfT_{\mcV})= 
\Loc\big(\{\kos C{\fp} \mid \text{$C\in \sfG$ and $\fp\in \mcV$}\}\big)\,.\]
  In particular, the subcategory $\sfT_{\mcV}$ is generated by a subset of $\sfT^c$.
\end{theorem} 

\begin{proof}
Let $\sfS=\Loc\big(\{\kos C{\fp} \mid \text{$C\in \sfG$ and $\fp\in \mcV$}\}\big)$. We then have inclusions
\[ 
\sfS \subseteq \Loc(\sfT^c\cap \sfT_{\mcV}) \subseteq \sfT_{\mcV}\,,
\] 
because $\kos C{\fp}$ belongs to $\sfT_{\mcV}$ for each $\fp$ in $\mcV$, by Lemma~\ref{le:koszul-properties}(2), and $\sfT_{\mcV}$ is a localizing subcategory, by Lemma~\ref{le:localizing}. It remains to prove that $\sfS=\sfT_{\mcV}$. Since the category $\sfS$ is compactly generated, the inclusion $\sfS\subseteq \sfT_{\mcV}$ admits a right adjoint $F\col \sfT_{\mcV} \to \sfS$.  Fix  $X\in \sfT_{\mcV}$ and complete the adjunction morphism $FX\to X$ to an exact triangle
\[
FX \lto X \lto Y\lto\,.
\]  
Given Corollary~\ref{cor:support}, the desired result follows once we verify that $\supp_{R}Y=\varnothing$. 
Fix $\fp\in \mcV$ and an object $C$ in $\sfG$. Since $\kos C{\fp}$ is in $\sfS$, the map $F{X}\to X$ induces an isomorphism
\[
\Hom_{\sfT}^*(\kos C{\fp},FX)\cong \Hom_{\sfT}^*(\kos C{\fp},X)\,.\]
  This implies $\Hom_{\sfT}^*(\kos C{\fp},Y)=0$, and therefore $\Hom_{\sfT}^*(C,\gam_{\fp}Y)=0$, by Proposition~\ref{pr:koszul-support}. Since this holds for each $C$ in $\sfG$, and $\sfG$ generates $\sfT$, one obtains that $\Hom^*_{\sfT}(-,\gam_{\fp}Y)=0$ on $\sfT$ and hence that $\gam_{\fp}Y=0$. Thus
\[
\supp_{R}Y\subseteq (\spec R)\setminus \mcV\,.\]
  Corollary~\ref{cor:support} yields $\supp_{R}Y\subseteq \mcV$, as $FX$ and $X$ are in $\sfT_{\mcV}$, so $\supp_{R}Y=\varnothing$.
\end{proof} 

The theorem above implies that local cohomology functors and localization functors are \emph{smashing}, that is to say, they preserves small coproducts: 

\begin{corollary}\label{cor:localization-smashing}
Let $\mcV\subseteq \spec R$ be specialization closed.  The exact functors $\gam_{\mcV}$ and $\bloc_{\mcV}$ on $\sfT$  preserve small coproducts.
\end{corollary}

\begin{proof}
  By the preceding theorem, $\sfT^c\cap\sfT_{\mcV}$ generates $\sfT_{\mcV}$, hence an  object $Y$ in $\sfT$ is $L_\mcV$-local if and only if $\Hom_\sfT(-,Y)=0$ on  $\sfT^c\cap\sfT_{\mcV}$. Therefore, the subcategory of $L_\mcV$-local objects is closed under taking small coproducts, and hence $L_\mcV$ preserves small coproducts, see Lemma~\ref{le:local-acyclic}. Using the localization triangle connecting $\gam_{\mcV}$ and $L_\mcV$, it follows that $\gam_{\mcV}$ preserves small coproducts as well.
\end{proof} 

\begin{corollary}\label{cor:support-sums}
\pushQED{\qed}%
Given a set of objects $X_i$ in $\sfT$, one has an equality
\[
\supp_{R}\coprod_i X_i=\bigcup_i\supp_{R} X_i\,. \qedhere\]
  \end{corollary} 

\subsection*{A recollement}
The functors $\gam_{\mcV}$ and $\bloc_{\mcV}$ corresponding to a specialization closed subset $\mcV\subseteq\spec R$ form part of a recollement.  Recall, see \cite[Sect.~1.4]{BBD},  that a \emph{recollement} is a diagram of exact functors
\[
\xymatrix{\sfT'\ar[rr]|-{I}&&\sfT \ar[rr]|-{Q}
\ar@<1.5ex>[ll]^-{I_\lambda}\ar@<-1.5ex>[ll]_-{I_\rho}&&\sfT''
\ar@<1.5ex>[ll]^-{Q_\lambda}\ar@<-1.5ex>[ll]_-{Q_\rho}}\] satisfying
the following conditions.
\begin{enumerate}
\item $I_\lambda$ is a left adjoint and  $I_\rho$ a right adjoint of $I$;
\item $Q_\lambda$ is a left adjoint and  $Q_\rho$ a right adjoint of $Q$;
\item $I_\lambda I\cong\Id_{\sfT'}\cong I_\rho I$ and $QQ_\rho \cong\Id_{\sfT''}\cong QQ_\lambda$;
\item $\Im I=\Ker Q$, that is, $QX=0$ holds iff $X\cong IX'$ for some $X'$ in $\sfT'$.
\end{enumerate}

Given any subset $\mcU\subseteq\spec R$, we denote by $\sfT(\mcU)$ the
full subcategory of $\sfT$ which is formed by all objects $X$ with
$\supp_{R} X\subseteq \mcU$.  Observe that this subcategory is not
necessarily the same as $\sfT_{\mcU}$, introduced in
Section~\ref{se:cohomology}; however, they coincide when $\mcU$ is
specialization closed, see Corollary~\ref{cor:support}.

\begin{theorem}\label{thm:recollement}
Let $\mcV\subseteq\spec R$ be specialization closed and set $\mcU=\spec R\setminus \mcV$.
The inclusion functor $u\colon\sfT(\mcU)\to\sfT$ then induces the following recollement
\[
\xymatrix{\sfT(\mcU)\ar[rr]|-{u=\inc}&&\sfT \ar[rr]|-{v=\gam_\mcV}
\ar@<1.5ex>[ll]^-{u_\lambda=\bloc_\mcV}\ar@<-1.5ex>[ll]_-{u_\rho}&&\sfT(\mcV)=\sfT_{\mcV}
\ar@<1.5ex>[ll]^-{v_\lambda=\inc}\ar@<-1.5ex>[ll]_-{v_\rho}}\]
with $u_\lambda$ a left adjoint of $u$, $u_\rho$ a right adjoint of $u$,
$v_\lambda$ a left adjoint of $v$, and $v_\rho$ a right adjoint of $v$.
\end{theorem} 

\begin{proof}
Since $\sfT(\mcU)$ equals the subcategory of $\bloc_\mcV$-local objects and $\sfT(\mcV)$ equals the subcategory of $\bloc_\mcV$-acyclic objects, by Lemma~\ref{le:local-acyclic}, the functors $\bloc_\mcV$ and $\gam_\mcV$ give rise to the following diagram:
\[
\xymatrix{\sfT(\mcU)\ar@<.75ex>[rr]^-{u=\inc}&&\sfT \ar@<.75ex>[rr]^-{v=\gam_\mcV}
\ar@<.75ex>[ll]^-{u_\lambda=\bloc_\mcV}&&\sfT(\mcV)=\sfT_{\mcV}
\ar@<.75ex>[ll]^-{v_\lambda=\inc}\,.}
\]
We have seen in Corollary~\ref{cor:localization-smashing} that the
functors $\bloc_\mcV$ and $\gam_\mcV$ preserves small coproducts. Thus
$u$ and $v$ preserve small coproducts. Then Brown
representability implies that $u$ and $v$ admit right adjoints
$u_\rho$ and $v_\rho$, respectively. It is straightforward to check
that these functors satisfy the defining conditions of a recollement.
\end{proof} 

\begin{remark}
In the notation of the preceding theorem, the functor
$v_\rho\comp\gam_{\mcV}$ on $\sfT$, which is right adjoint to
$\gam_{\mcV}$, may be viewed as a completion along $\mcV$, see
Remark~\ref{rem:completions}.
\end{remark} 

\section{Connectedness}
\label{se:connectedness}
As before, let $\sfT$ be a compactly generated $R$-linear triangulated category. In this section we establish Krull-Remak-Schmidt type results for objects in $\sfT$, and deduce as a corollary the connectedness of supports of indecomposable objects. We give a second proof of this latter result, by deriving it from an analogue of the classical
Mayer-Vietoris sequence in topology. 

\begin{theorem}\label{thm:KRS}
If $X\in\sfT$ is such that $\supp_{R}X\subseteq \bigsqcup_{i\in I} \mcV_i$ where the subsets $\mcV_i$ are pairwise disjoint and specialization closed, then there is a natural isomorphism:
\[
X\cong \coprod_{i\in I}\gam_{\mcV_i}X\,.\]
  \end{theorem} 

\begin{proof}
  The canonical morphisms $\gam_{\mcV_i}X\to X$ induce a morphism $\eps$ as below
\[
\coprod_{i\in I} \gam_{\mcV_i}X \xra{\  \eps \ } X \lto Y\lto\]
  We then complete it to an exact triangle as above. One has $\supp_{R}Y\subseteq \supp_{R}X$; this follows from Theorem~\ref{thm:localization-support}, Proposition~\ref{pr:triangle}, and Corollary~\ref{cor:support-sums}. 

We claim that $\supp_{R}Y \cap \supp_{R}X=\varnothing$, which then implies $Y=0$, by Theorem~\ref{thm:detection}, and hence that $\eps$ is an isomorphism, as desired. 

Indeed, pick a point $\fp$ in $\supp_{R}X$. There is a unique index $k$ in $I$ for which $\fp$ is in $\mcV_k$. Applying $\gam_{\fp}(-)$, and keeping in mind Corollary~\ref{cor:localization-smashing}, one obtains a diagram
\[
\gam_{\fp}\gam_{\mcV_k}X \xra{\ \cong \ }
      \coprod_{i\in I} \gam_{\fp}\gam_{\mcV_i}X  \xra{\ \gam_{\fp}\eps\ } \gam_{\fp}X\,,\]
  The isomorphism holds because $\gam_{\fp}\gam_{\mcV_i}X=0$ for $i\ne k$, by Proposition~\ref{pr:localizations-commute}. The same result yields also that the composed map $\gam_\fp\gam_{\mcV_k}X\to \gam_\fp X$ is an isomorphism, so one deduces that $\gam_{\fp}\eps$ is also an isomorphism. Therefore, the exact triangle above implies $\gam_{\fp}Y=0$, that is to say, $\fp\notin\supp_{R}Y$.
\end{proof} 

Recall that a specialization subset $\mcV$ of $\spec R$ is said to be \emph{connected} if for any pair
$\mcV_{1}$ and $\mcV_{2}$ of specialization closed subsets of $\spec R$, one has
\[
\mcV\subseteq\mcV_{1}\cup\mcV_{2}\text{ and } \mcV_{1}\cap\mcV_{2}=\varnothing\quad
\Longrightarrow\quad \mcV\subseteq\mcV_{1}\text{ or }\mcV\subseteq\mcV_{2}\,.\]
  
The following lemma is easy to prove.
   
\begin{lemma}\label{le:KRS}
Each specialization closed subset $\mcV$ of $\spec R$ admits a unique decomposition $\mcV= \bigsqcup_{i\in I}
\mcV_i$ into nonempty, specialization closed, connected, and pairwise disjoint subsets. \qed
\end{lemma} 

Theorem~\ref{thm:KRS} yields Theorem~\ref{intro:KRS} stated in the introduction. 

\begin{theorem}\label{thm:KRS2}
Each object $X$ in $\sfT$ admits a unique decomposition $X=\coprod_{i\in I} X_i$ with $X_i\ne 0$ such that the subsets $\cl(\supp_{R} X_i)$ are connected and pairwise disjoint.
\end{theorem} 

\begin{proof}
Use Lemma~\ref{le:KRS} to get a decomposition $\cl(\supp_{R}X)=\bigsqcup_{i\in I} \mcV_i$ into connected, pairwise disjoint specialization closed subsets, and then apply Theorem~\ref{thm:KRS} to obtain a decomposition $X=\coprod_{i\in I} X_i$, where $X_i=\gam_{\mcV_i}X$.  Observe that $X_i\ne 0$ by Corollary~\ref{cor:support}, since $\mcV_i\cap \supp_{R}X \ne \varnothing$.  It is easy to verify the other properties of
the decomposition. 
  
Let $X=\coprod_{j\in J} Y_j$ be another such decomposition.  Lemma~\ref{le:KRS} then implies that there is a bijection $\sigma\colon I\to J$ with $\supp_{R} Y_{\sigma(i)}=\supp_{R} X_i$ for all $i$.  Corollary~\ref{cor:orthogonality} yields $\Hom_{\sfT}(X_i,Y_j)=0$ for $j\ne\sigma(i)$. Therefore, $X_i=Y_{\sigma(i)}$ for each $i$.
\end{proof} 

The connectedness theorem is a direct consequence of the preceding result. 

\begin{corollary}\label{cor:connectedness}
If $X$ is indecomposable, then $\cl(\supp_{R} X)$ is connected.\qed
\end{corollary} 

Following Rickard \cite{Ri}, one could deduce this result also from a simple special case
of a Mayer-Vietors triangle, described below. 

\subsection*{Mayer-Vietoris triangles}
When $\mcV\subseteq \mcW$ are specialization closed subsets of $\spec R$, one has natural morphisms:
\[
\gam_{\mcV}\xra{\ \theta_{\mcV,\mcW}\ } \gam_{\mcW}
\quad\text{and}\quad
\bloc_{\mcV}\xra{\ \eta_{\mcV,\mcW}\ }\bloc_{\mcW}\,.
\] 
We refer to the exact triangles in the next result as the \emph{Mayer-Vietoris triangles}
associated to $\mcV$ and $\mcW$. 

\begin{theorem}\label{thm:mayer-vietoris}
If $\mcV$ and $\mcW$ are specialization closed subsets of $\spec R$, then for each $X$ in $\sfT$, there are natural exact triangles:

\[  
\xymatrixcolsep{1.8pc}
\xymatrix{
\gam_{\mcV\cap\mcW}X \ar@{->}[rrr]^-{(\theta_{\mcV\cap\mcW,\mcV},\, \theta_{\mcV\cap\mcW,\mcW})^t}
  &&&\gam_{\mcV}X\amalg\gam_{\mcW}X \ar@{->}[rrr]^-{(\theta_{\mcV,\mcV\cup\mcW},\, -\theta_{\mcW,\mcV\cup\mcW})}
  &&&\gam_{\mcV\cup\mcW}X \ar@{->}[r]&\\
  \bloc_{\mcV\cap\mcW}X \ar@{->}[rrr]^-{(\eta_{\mcV\cap\mcW,\mcV},\, \eta_{\mcV\cap\mcW,\mcW})^t}
   &&&\bloc_{\mcV}X \amalg \bloc_{\mcW}X 
     \ar@{->}[rrr]^-{(\eta_{\mcV,\mcV\cup\mcW},\, -\eta_{\mcW,\mcV\cup\mcW})}
    &&& \bloc_{\mcV\cup\mcW}X  \ar@{->}[r]&{} }
\]
\end{theorem}

\begin{proof}
Complete the morphism
\[  
 \gam_\mcV(X)\amalg\gam_\mcW(X)
\xra{(\theta_{\mcV,\mcV\cup\mcW},-\theta_{\mcW,\mcV\cup\mcW})}\gam_{\mcV\cup\mcW}X
\]
to an exact triangle
\[  
  X'\lto \gam_\mcV(X)\amalg\gam_\mcW(X)
\xra{(\theta_{\mcV,\mcV\cup\mcW},-\theta_{\mcW,\mcV\cup\mcW})}\gam_{\mcV\cup\mcW}X\lto\,.
\]
It is easy to check that the morphism
\[  
  \gam_{\mcV\cap\mcW}X\xra{(\theta_{\mcV\cap\mcW,\mcV},\theta_{\mcV\cap\mcW,\mcW})^t}
   \gam_\mcV(X)\amalg\gam_\mcW(X)
\]
factors through $X'\to \gam_\mcV(X)\amalg\gam_\mcW(X)$, and gives a new exact triangle
\[  
  \gam_{\mcV\cap\mcW}X\lto X'\lto X''\lto\,.
\]
The goal is to prove that $X''=0$, so $\gam_{\mcV\cap\mcW}X\to X'$ is an isomorphism.  

First observe that $\supp_{R} X'$ and $\supp_{R} X''$ are contained in $\mcV\cup\mcW$, by Proposition~\ref{pr:triangle}.  Consider now the triangle below. We do not specify the third map and it is not asserted
that the triangle is exact.
\[  
  \gam_{\mcV\cap\mcW}X\xra{(\theta_{\mcV\cap\mcW,\mcV},\theta_{\mcV\cap\mcW,\mcW})^t}
\gam_\mcV(X)\amalg\gam_\mcW(X)
\xra{(\theta_{\mcV,\mcV\cup\mcW},-\theta_{\mcW,\mcV\cup\mcW})}\gam_{\mcV\cup\mcW}X\lto\,.
\]
It is readily verified that applying $\gam_\mcV$ to it yields a split exact triangle, and hence that $\gam_{\mcV\cap\mcW}X\to \gam_\mcV{X'}$ is an isomorphism. In particular, $\gam_\mcV{X''}=0$; in the same vein one deduces that $\gam_\mcW{X''}=0$.  Thus Corollary~\ref{cor:support} yields
\[  
  (\mcV\cup\mcW) \cap \supp_{R}{X''}=\varnothing\,.
\]
We conclude that $\supp_{R}{X''}=\varnothing$, and therefore $X''=0$ by Theorem~\ref{thm:detection}.  This establishes the first Mayer-Vietoris triangle. The proof for the second is similar.
\end{proof} 

\begin{proof}[Second proof of Corollary~\emph{\ref{cor:connectedness}}]
Suppose $\mcV$ and $\mcW$ are specialization closed subsets of $\spec R$ such that $\supp_{R} X\subseteq\mcV\cup\mcW$ and $\mcV\cap\mcW=\varnothing$. Then, since $\gam_{\varnothing}X=0$, the first Mayer-Vietoris triangle in Theorem~\ref{thm:mayer-vietoris} yields an isomorphism
\[  
  X\cong \gam_{\mcV\cup\mcW}{X}  \cong \gam_{\mcV}X\, \amalg \, \gam_{\mcW}X\,.
\]
In view of Theorem~\ref{thm:detection}, this implies that when $X$ is indecomposable, one of the subsets $\mcV\cap \supp_{R} X$ or $\mcW\cap \supp_{R}X$ is empty.
\end{proof} 

\section{Tensor triangulated categories}
\label{se:tensor}

Let $\sfT$ be a compactly generated triangulated category. In this section, $\sfT$ is also tensor
triangulated. Thus, $\sfT=(\sfT,\otimes,\bbi)$ is a symmetric monoidal
category with a \emph{tensor product}
$\otimes\col\sfT\times\sfT\to\sfT$ and a \emph{unit} $\bbi$. In
addition, we assume that the tensor product is exact in each variable
and that it preserves small coproducts; for details see
\cite[III.1]{LMS}. The Brown representability theorem yields
\emph{function objects} $\funct(X,Y)$ satisfying
\[  
  \Hom_\sfT(X\otimes Y,Z)\cong\Hom_\sfT(X,\funct(Y,Z))\quad\text{for all $X,Y,Z$ in $\sfT$}\,.
\]
For each $X$ in $\sfT$ we write
\[  
  X^\vee=\funct(X,\bbi)\,.
\]
It is assumed that the unit $\bbi$ is compact and that all compact objects $C$ are \emph{strongly dualizable}, that is, the canonical morphism
\[  
  C^\vee\otimes X\to \funct(C,X)
\]
is an isomorphism for all $X$ in $\sfT$.  

Let $C,D$ be compact objects in $\sfT$. The following properties are easily verified.
\begin{enumerate}[{\quad\rm(1)}]
\item $\bbi^\vee\cong\bbi$ and $C^{\vee\vee}\cong C$;
\item $\Hom_\sfT(X\otimes C^\vee,Y)\cong \Hom_\sfT(X,C\otimes Y)$, for all $X,Y$ in $\sfT$;
\item $C^\vee$ and $C\otimes D$ are compact.\end{enumerate}These properties are used in the sequel without further comment. 

\subsection*{Tensor ideals and smashing localization}
A full subcategory $\sfS$ of $\sfT$ is called a \emph{tensor ideal} in $\sfT$ if for each $X$ in $\sfS$ and $Y$ in $\sfT$, the object $X\otimes Y$ belongs to $\sfS$. This condition is equivalent to: $Y\otimes X$ belongs to $\sfS$, as the tensor product is symmetric. 

\begin{proposition}\label{pr:smash-loc}
Let $\bloc\col\sfT\to\sfT$ be a localization functor such that the category $\sfT_{L}$ of $\bloc$-acyclic objects is generated by $\sfT_{L}\cap\sfT^c$ and the latter is a tensor ideal in $\sfT^c$. Then the following statements hold. 
\begin{enumerate}[{\quad\rm(1)}]
\item The $\bloc$-acyclic objects and the $\bloc$-local objects are both tensor ideals in $\sfT$.
\item For each $X$ in $\sfT$, one has natural isomorphisms
\[  
  \gam{X}\cong X\otimes\gam{\bbi}\quad\text{and}\quad\bloc{X}\cong X\otimes \bloc{\bbi}\,.
\]
\end{enumerate}
\end{proposition} 

\begin{proof}
    (1) Let $X$ be an object in $\sfT$. One  has the following equivalences:
\begin{align*}
X\text{ is $\bloc$-local} &\iff \Hom_\sfT(C,X)=0\text{ for all }C\in\sfT_{L}\cap\sfT^c\\ 
&\iff \Hom_\sfT(C\otimes D^\vee,X)=0\text{ for all }C\in\sfT_{L}\cap\sfT^c 
    \text{ and }D\in\sfT^c\\ 
&\iff \Hom_\sfT(C,X\otimes D)=0\text{ for all }C\in\sfT_{L}\cap\sfT^c \text{ and }D\in\sfT^c\\ &\iff \Hom_\sfT(C,X\otimes Y)=0\text{ for all }C\in\sfT_{L}\cap\sfT^c\text{ and }Y\in\sfT\\ 
&\iff X\otimes Y\text{ is $\bloc$-local for all }Y\in\sfT\,.
\end{align*}
The first and the last equivalences hold because $\sfT_{L}\cap\sfT^c$ generates $\sfT_{L}$.  The second holds because $\bbi\cong \bbi^\vee$ is compact, and $C\otimes D^{\vee}$ is compact and $L$-acyclic; the third one holds because $D$ is strongly dualizable; the fourth holds because the functor $\Hom_\sfT(C,X\otimes -)$ preserves exact triangles and small coproducts, and $\sfT$ is compactly generated.  Thus, the $\bloc$-local objects form a tensor ideal in $\sfT$. 

Consider now the $\bloc$-acyclic objects.  As the $\bloc$-local objects form a tensor
ideal in $\sfT$, one obtains the second step in the following chain of equivalences:
\begin{align*}
X\text{ is $\bloc$-acyclic} &\iff \Hom_\sfT(X,Z)=0\text{ for all $\bloc$-local }Z\in\sfT\\ 
&\iff \Hom_\sfT(X,Z\otimes C^\vee)=0\text{ for all }C\in\sfT^c\text{ and $\bloc$-local }Z\in\sfT\\
&\iff \Hom_\sfT(X\otimes C,Z)=0\text{ for all }C\in\sfT^c\text{ and $\bloc$-local }Z\in\sfT\\
&\iff \Hom_\sfT(X\otimes Y,Z)=0\text{ for all }Y\in\sfT\text{ and $\bloc$-local }Z\in\sfT\\
&\iff X\otimes Y\text{ is $\bloc$-acyclic for all }Y\in\sfT\,.
\end{align*}
The justifications of the other equivalences is similar to those in the previous
paragraph. Therefore, the $\bloc$-acyclic objects form a tensor ideal. 

(2) Let $\eta\col\Id_{\sfT}\to \bloc$ be the morphism associated to the localization
functor $\bloc$. Let $X$ be an object in $\sfT$, and consider the commutative square
\[  
  \xymatrix{
X\ar[rr]^-{X\otimes\eta{\bbi}}\ar[d]^-{\eta{X}} &&
          X\otimes\bloc{\bbi}\ar[d]^-{\eta(X\otimes \bloc{\bbi})}\\ 
     \bloc{X}\ar[rr]^-{\bloc(X\otimes\eta{\bbi})}&& 
          \bloc(X\otimes \bloc{\bbi})}
\]
One has $\bloc(X\otimes\gam{\bbi})=0$ and $\gam(X\otimes\bloc{\bbi})=0$, because the
$\bloc$-acyclic and the $\bloc$-local objects form tensor ideals, by part (1). Therefore
$\bloc(X\otimes\eta{\bbi})$ and $\eta(X\otimes \bloc{\bbi})$ are isomorphisms, and hence
$\bloc{X}\cong X\otimes\bloc{\bbi}$. 

A similar argument shows that $\gam{X}\cong X\otimes\gam{\bbi}$.
\end{proof}  

\subsection*{Support}
In the remainder of this section we fix, as before, a homomorphism of graded rings $R\to Z(\sfT)$ into the graded center of $\sfT$.  Note that the endomorphism ring $\End^*_\sfT(\bbi)$ is graded-commutative, and that the homomorphism
\[  
  \End^*_\sfT(\bbi)\lto Z(\sfT), \quad\text{where}\quad\alpha\mapsto\alpha\otimes -\,,
\]
is a canonical choice for the homomorphism $R\to Z(\sfT)$. 

\begin{theorem}
  Let $\mcV$ be a specialization closed subset of $\spec R$. Then the $\bloc_\mcV$-acyclic
  objects and the $\bloc_\mcV$-local objects both form tensor ideals in $\sfT$. Moreover,
  we have for each $X$ in $\sfT$ natural isomorphisms
\[  
  \gam_{\mcV}X\cong X\otimes\gam_\mcV{\bbi} \quad\text{and}\quad
  \bloc_{\mcV}X\cong X\otimes \bloc_\mcV{\bbi}\,.
\]
\end{theorem} 

\begin{proof}
    We know from Theorem~\ref{thm:localization-finite} that the subcategory $\sfT_{\mcV}$ of
  $\bloc_\mcV$-acyclic objects is generated by compact objects. Thus the assertion follows
  from Proposition~\ref{pr:smash-loc}, once we check that $\sfT_{\mcV}\cap\sfT^c$ is a
  tensor ideal in $\sfT^c$. This last step is contained in the following chain of
  equivalences:
\begin{align*}
X\in\sfT_{\mcV} &\iff\supp_R\Hom^*_\sfT(B,X)\subseteq\mcV\text{ for all }B\in\sfT^c\\
&\iff\supp_R\Hom^*_\sfT(B\otimes C^\vee,X)\subseteq\mcV\text{ for all }B,C\in\sfT^c\\
&\iff\supp_R\Hom^*_\sfT(B,X \otimes C)\subseteq\mcV\text{ for all }B,C\in\sfT^c\\
&\iff X\otimes C\in\sfT_{\mcV} \text{ for all }C\in\sfT^c\,.\qedhere
\end{align*}\qedhere
\end{proof} 

The following corollary is immediate from the definition of $\gam_{\fp}$: 

\begin{corollary}
\pushQED{\qed}%
Let $\fp$ be a point in $\spec R$. For each object $X$ in $\sfT$ one has a natural isomorphism
\[  
  \gam_{\fp}X\cong X\otimes\gam_{\fp}{\bbi}\,.\qedhere
\]
\end{corollary} 

\begin{remark}
  Setting $R=\End_\sfT^*(\bbi)$, one recovers the notion of support for noetherian stable
  homotopy categories discussed in \cite[Sect.~6]{HPS}.  Theorem~\ref{thm:KRS2}
  specialized to this context extends the main result in \cite{Ba}; see also
  \cite{Ca,Kr:krs}.
\end{remark} 

For commutative noetherian rings, Hopkins~\cite{Ho} and
Neeman~\cite{Ne:dr} have classified thick subcategories of perfect
complexes, and localizing subcategories of the derived category, in
terms of subsets of the spectrum.  These have been extended to the
realm of algebraic geometry by Thomason~\cite{Th}, and Alonso Tarrío,
Jerem\'ias L\'opez, and Souto Salorio~\cite{JVS2}. In \cite{HPS}, some
of these results have been extended to tensor triangulated
categories. However, this context do not cover important examples;
notably, the bounded derived category of a complete intersection ring;
see Section~\ref{se:ci}.  One of the motivations for this article was
to develop a broader framework, and attendant techniques, where we
could state and prove such results.

\section{Commutative noetherian rings}
\label{se:commutative}
In this section we apply the theory developed in the preceding
sections to the derived category of a commutative noetherian ring. The
main result interprets localization functors with respect to
specialization closed subsets to the corresponding local cohomology
functors, and establishes that, in this context, the notion of support
introduced here coincides with the classical one.

Throughout this section, $A$ is a commutative noetherian ring and
$\sfM=\Mod A$ the category of $A$-modules.  Let $\sfT=\sfD(\sfM)$ be
the derived category of complexes of $A$-modules.  The category $\sfT$
is triangulated, and admits coproducts.  The module $A$, viewed as a
complex concentrated in degree $0$, is a compact generator for $\sfT$.
In what follows, for each complex $X$ of $A$-modules $H^*(X)$ denotes
the cohomology of $X$. Note that
\[  
  H^*(X)=\Hom^*_{\sfT}(A,X)\,,
\]
so the notation $H^*(-)$ is compatible with its usage in Theorem~\ref{thm:axioms}. 
Next we recall the notion of local cohomology and refer to \cite{Ha,Li} for details. 

\subsection*{Local cohomology}
Let $\mcV$ be a specialization closed subset of $\spec A$, and $M$ an $A$-module.
Consider the submodule $F_{\mcV}M$ of $M$ defined by the exact sequence:
\[  
  0\lto F_{\mcV}M\lto M \lto \prod_{\fq\notin\mcV}M_\fq\,.
\]
The traditional notation for $F_{\mcV}M$ is $\gam_{\mcV}M$, but we
reserve that for the derived version, see Theorem~\ref{thm:ca-main}.
The assignment $M\mapsto F_{\mcV}M$ is an additive, left-exact functor
on the category of $A$-modules, and inclusion gives a morphism of
functors $F_{\mcV}\to \Id_\sfM$.  Observe that $F_\mcV$ provides a
right adjoint to the inclusion $\sfM_\mcV\subseteq\sfM$; this follows from
Lemma~\ref{le:specialization}. 

The following properties of $F_{\mcV}$ are readily verified using the
exact sequence above.
\begin{enumerate}[{\quad\rm(1)}]
\item For any ideal $\fa$ in $R$, the closed subset $\mcZ=\mcV(\fa)$ satisfies
\[  
  F_{\mcZ}M=\{m\in M\mid \text{$\fa^n m=0$ for some integer $n\ge 0$} \}\,.
\]
\item For each arbitrary specialization closed subset $\mcV$, one has that
\[  
  F_{\mcV}M  = \bigcup_{\substack{\mcZ\subseteq \mcV\\ \text{closed}}}F_{\mcZ}M \,.
\]
\item For each prime ideal $\fp$ in $A$ one has that
\[  
  F_{\mcV}(E(A/\fp)) = \begin{cases}
E(A/\fp) & \text{if   $\fp\in\mcV$}\\
0        & \text{othewise.}     
\end{cases}
\]  
\end{enumerate} 

We denote by $\lch {\mcV}\col \sfT\to \sfT$ the right derived functor
of $F_{\mcV}$. For each complex $X$ of $A$-modules the \emph{local
  cohomology} of $X$ with respect to $\mcV$ is the graded $A$-module
\[  
  \coh{\mcV}X = \coh{}{\lch{\mcV}X}\,.
\]
Thus, if $I$ is an injective resolution of $X$, then
$\coh{\mcV}X=\coh{}{F_{\mcV}I}$, where $F_{\mcV}I$ is the complex of
$A$-modules with $(F_{\mcV}I)^n=F_{\mcV}(I^n)$, and differential
induced by the one on $I$.  When $\fa$ is an ideal in $R$, it is
customary to write $\coh{\fa}X$ for the local cohomology of $X$ with
respect to the closed set $\mcV(\fa)$.  We consider support with
respect to the canonical morphism
\[  
  A\lto Z(\sfT)
\]
given by homothety: $a\mapsto a\cdot \id$; here $A$ is viewed as a
graded ring concentrated in degree zero.  The next theorem explains
the title of this paper.  In the sequel, $X_\fp$ denotes the complex
of $A_\fp$-modules $A_\fp\otimes_AX$.

\begin{theorem}
\label{thm:ca-main}
For each specialization closed subset $\mcV$ of $\spec A$, one has an isomorphism
$\gam_{\mcV}\cong \lch{\mcV}$.  Moreover, each complex $X$ of $A$-modules satisfies
\[  
  \supp_{A}X = \{\fp\in \spec A\mid \coh {\fp A_\fp}{X_\fp}\ne 0\}\,.
\]
\end{theorem} 

\begin{remark}
  Foxby~\cite{Fo} has proved that for each complex $X$ of $A$-modules
  with $H^*(X)$ bounded, and each $\fp\in\spec A$, the following
  conditions are equivalent:
\begin{enumerate}[{\quad\rm(1)}]
\item $E(A/\fp)$ occurs in the minimal injective resolution of $X$;
\item $\Ext_{A_\fp}^*(k(\fp),X_\fp)\neq 0$;
\item $\Tor^A_{*}(k(\fp),X)\neq 0$;
\item $\hH_{\fp A_\fp}^{*}(X_\fp)\ne 0$.
\end{enumerate}
Therefore, the theorem above implies that for modules, the notion of
support as defined in this section, coincides with the one from
Section~\ref{se:modules}.  Moreover, it is implicit in the results in
\cite{Ne:dr}, and is immediate from \cite[(2.1) and (4.1)]{FI}, that
the conditions (1)--(4) coincide for every complex $X$ of modules,
which means that the set $\supp_AX$ is a familiar one.
\end{remark} 

\begin{proof}[Proof of Theorem~\emph{\ref{thm:ca-main}}]
We let $\sfM_\mcV$ denote the full subcategory of $\sfM$ consisting of
modules with support, in the sense of Section~\ref{se:modules},
contained in $\mcV$.  In particular,
$\sfT_\mcV=\sfD_{\sfM_\mcV}(\sfM)$, the full subcategory of $\sfT$
formed by complexes $X$ with $H^n(X)$ in $\sfM_\mcV$, for each
$n$. The inclusion functor $\sfM_\mcV\to\sfM$ induces an equivalence
of categories $\sfD(\sfM_\mcV)\to \sfD_{\sfM_\mcV}(\sfM)$, because
every module in $\sfM_\mcV$ admits a monomorphism into a module in
$\sfM_\mcV$ which is injective in $\sfM$; see
\cite[Proposition~I.4.8]{Ha}.  

Next observe that the functor $\sfD(\sfM_\mcV)\to\sfD(\sfM)$ equals
the left derived functor of the inclusion $\sfM_\mcV\to\sfM$. Thus the
right derived functor $\lch{\mcV}$ provides a right adjoint to the
functor $\sfD(\sfM_\mcV)\to\sfD(\sfM)$, since $F_\mcV$ is a right
adjoint to the inclusion $\sfM_\mcV\to\sfM$. On the other hand,
by construction $\gam_{\mcV}$ is a right adjoint to the inclusion
$\sfD_{\sfM_\mcV}(\sfM)\to\sfD(\sfM)$.  A right adjoint functor is
unique up to isomorphism, so the equivalence $\sfD(\sfM_\mcV)\cong
\sfD_{\sfM_\mcV}(\sfM)$ implies $\lch{\mcV}\cong\gam_{\mcV}$.

Let $\fp$ be a point in $\spec A$, and set $\mcZ(\fp)=\{\fq\in\spec
A\mid \fq\not\subseteq\fp\}$. Let $X$ be a complex of $A$-modules.  The
functor $\sfT\to\sfT$ sending $X$ to $X_\fp$ is a localization functor
and has the same acyclic objects as $\bloc_{\mcZ(\fp)}$, since
$H^*(\bloc_{\mcZ(\fp)}X)\cong H^*(X)_\fp\cong H^*(X_\fp)$.  For the
first isomorphism, see Theorem~\ref{thm:locprime}. Therefore
$\bloc_{\mcZ(\fp)}X\cong X_\fp$, so one has isomorphisms
\[  
  \gam_{\fp}X = \gam_{\mcV(\fp)}\bloc_{\mcZ(\fp)}X 
              \cong \gam_{\mcV(\fp)}X_\fp 
              \cong \lch{\mcV(\fp)}X_\fp\,.
\]
This yields the stated expression for $\supp_{R}X$, as
$H^*(\lch{\mcV(\fp)}X_\fp)=\coh{\fp A_\fp}{X_\fp}$.
\end{proof} 

\begin{remark}
  One can give other proofs for the first part of the preceding
  theorem. For instance, it is easy to check the functor $F_{\mcV}$
  coincides with $\gam_{\mcV}'$ in \cite[(3.5)]{Li}, so Proposition
  3.5.4 in \emph{loc.~cit.} yields that $\lch{\mcV}$ is right adjoint
  to the inclusion $\sfT_{\mcV}\subseteq\sfT$. It must thus coincide
  with $\gam_{\mcV}$. One can also approach this result via the
  machinery in Section~\ref{se:tensor}, for the derived category is
  tensor triangulated.

  Here is a different perspective: As before, one argues that $\bloc_{\mcZ(\fp)}X\cong
  X_\fp$ for each complex $X$.  We claim that for any $A$-module $M$, if the support of
  $M$, as computed from the minimal injective resolution, equals $\{\fp\}$, then
  $\supp_{A}M = \{\fp\}$. 

  Indeed, $M\in \sfT_{\mcV(\fp)}$, by definition, so $\supp_{R}M
  \subseteq \mcV(\fp)$.  Fix a prime ideal $\fq\supset \fp$. The
  hypothesis on $M$ implies $M\cong M_\fp$, so one obtains the first
  isomorphism below.
\[  
  \gam_{\fq}M = \gam_{\mcV(\fq)}\bloc_{\mcZ(\fq)} \bloc_{\mcZ(\fp)}M
 \cong\gam_{\mcV(\fq)}\bloc_{\mcZ(\fp)}M=0\,.
\]
The remaining isomorphisms hold by Lemma~\ref{le:functors-commute},
since $\mcZ(\fp)$ contains $\mcV(\fq)$ and $\mcZ(\fp)$. Therefore,
$\fq\notin\supp_{A}M$. This settles the claim.

Let $\mcV$ be a specialization closed subset of $\spec A$.  Observe
that both $\lch{\mcV}$ and $\gam_{\mcV}$ are exact functors on $\sfT$,
and that $\lch{\mcV}X\in \sfT_{\mcV}$, as $\supp_A
H^*(\lch{\mcV}X)\subseteq \mcV$.  Since $\gam_{\mcV}$ is right adjoint
to the inclusion $\sfT_\mcV\subseteq \sfT$, one has thus a morphism
$\lch{\mcV}\to \gam_{\mcV}$. This is an isomorphism on injective
modules $E(A/\fp)$, for each $\fp\in\spec A$, by the preceding claim
and the properties of $F_{\mcV}$ listed above. It follows that
$\lch{\mcV}\to \gam_{\mcV}$ is an isomorphism on all of $\sfT$,
because $\lch{\mcV}$ and $\gam_{\mcV}$ commute with coproducts and
$\sfT$ equals the localizing subcategory generated by the $E(A/\fp)$,
see \cite[\S 2]{Ne:dr}.
\end{remark} 

The example below is intended to illustrate the difference between
$\supp_AX$ and $\supp_AH^*(X)$. In particular, we see that the
inclusions in Corollary~\ref{cor:support-inclusions} can be strict.
One can construct such examples over any commutative noetherian ring
of Krull dimension at least two; see \cite{Kr:thick}. 

\begin{example}
\label{ex:supp-incl}
Let $k$ be a field, let $A=k[\![x,y]\!]$, the power series ring in
indeterminates $x,y$, and set $\fm=(x,y)$, the maximal ideal of $A$.
The minimal injective resolution of $A$ has the form:
\[  
  \cdots\lto 0\lto Q\lto \bigoplus_{\height\fp=1} E(A/\fp)\lto E(A/\fm)\lto 0\lto \cdots
\]
where $Q$ denotes the fraction field of $A$. Let $X$ denote the truncated complex
\[  
  \cdots\lto 0\lto Q\lto \bigoplus_{\height\fp =1}E(A/\fp)\lto 0\lto\cdots
\]
viewed as an object in $\sfT$, the derived category of $A$. One then has
\[  
  \supp_{A} X=(\spec A)\setminus\{\fm\}\,, \quad
\Min_A H^*(X) = \{(0)\}\,,\quad  \supp_A H^*(X) =\spec A\,.
\]
Compare this calculation with the conclusion of Corollary~\ref{cor:support-inclusions}.
\end{example} 

In view of Theorems~\ref{thm:ca-main}, results on support established
in previous sections specialize to the case of complexes over
commutative noetherian rings. Two of these are noteworthy and are
commented upon.

\begin{remark}
  Corollary~\ref{cor:connectedness} yields that the support of an
  indecomposable complex of $A$-modules is connected. Restricting the
  decomposition in Theorem~\ref{thm:KRS2} to compact objects gives a
  Krull-Remak-Schmidt type theorem for thick subcategories of the
  category of perfect complexes.  This recovers results of
  Chebolu~\cite[(4.13), (4.14)]{Ch}.
\end{remark} 

\begin{remark}
\label{rem:completions}
Let $\fa$ be an ideal in $A$ and set $\mcV=\mcV(\fa)$. Consider
Theorem~\ref{thm:recollement} in the context of this section. The
functor $v_\rho\gam_{\mcV}$ is a right adjoint to $v_{\lambda}\gam_{\mcV}$. 
Given Theorem~\ref{thm:ca-main}, it follows from Greenlees-May duality---see \cite[??]{GM} and also \cite[\S 4]{Li}---that $v_\rho\gam_{\mcV}$ is the left derived of the $\fa$-dic completion
functor.
\end{remark} 

We have focussed on the derived category of $A$-modules. However,
similar considerations apply also to the homotopy category of
injective modules, and to the homotopy category of projective
modules. This leads to a notion of support for acyclic and totally
acyclic complexes in either category. This has connections to results
in \cite{IK}. We intend to pursue this line of investigation
elsewhere.

\section{Modules over finite groups}
\label{se:groups}

Let $G$ be a finite group and $k$ a field of characteristic $p$
dividing $|G|$. There are several choices of a tensor
triangulated category $\sfT$, and we comment on each of them. In each
case, we take for the ring $R$ the cohomology ring $H^*(G,k)$; by the
Evens--Venkov theorem~\cite[(3.10)]{Benson:1991} this is a finitely
generated graded-commutative $k$-algebra, and hence noetherian. We use
the tensor product $\otimes_k$ for objects in $\sfT$ with the usual
diagonal $G$-action: $g(x\otimes y)=gx\otimes gy$ for all $g\in G$
and $x\otimes y\in X\otimes_k Y$.

\subsection*{The stable module category}
The first choice for $\sfT$ is the stable module category $\StMod
kG$. In this case, the trivial module $k$ is compact, and plays the
role of the unit $\bbi$ in Section \ref{se:tensor}.  Note that $k$
does not necessarily generate $\sfT$. The function objects are
$\Hom_k(M,N)$ with the usual $G$-action: $(g(\alpha))(m)=g(\alpha(g^{-1}m))$.
The compact objects are the modules isomorphic in $\StMod kG$ to
finitely generated modules, and the full subcategory of compact
objects in $\StMod kG$ is denoted $\stmod kG$. Compact objects are
strongly dualizable.  The graded endomorphism ring of $k$ is
isomorphic to the Tate cohomology ring $\hat H^*(G,k)$.  So this maps
to the center of $\StMod kG$, but it is usually not
noetherian. Furthermore, the full center of $\StMod kG$ appears not to
be well understood in general.

\begin{lemma}
Let $G$ be a finite group and $k$ a field of characteristic $p$.
Then the following are equivalent:
\begin{enumerate}[\quad\rm(1)]

\item

$\hat H^*(G,k)$ is noetherian.
\item
$\hat H^*(G,k)$ is periodic.
\item
$G$ has $p$-rank one.
\end{enumerate}
\end{lemma}

\begin{proof}
If $G$ has $p$-rank one then $\hat H^*(G,k)$ is periodic and
noetherian, by \cite[(XII.11)]{CE}.  If $G$ has $p$-rank greater than
one then the negative Tate cohomology is contained in the nil radical
by Benson and Krause \cite[(2.4)]{BK:2002a}\footnote{Let us set the
record straight at this stage about the proof of this
proposition. Lemma 2.3 of that paper as stated and its proof are
obviously incorrect.  The correct statement, which is the one used in
the proof of 2.4, is that if $G$ has $p$-rank greater than one, and
$H$ is a subgroup of $p$-rank one, then the restriction to $H$ of any
element of \emph{negative} degree in $\hat H^*(G,k)$ is nilpotent. To
see this, if $x$ is such an element whose restriction is not
nilpotent, and $y\in \hat H^*(H,k)$ satisfies $\mathop{\rm
res}_{G,H}(x)y=1$, then $y$ has positive degree, so using the Evens
norm map, some power of $y$ is in the image of restriction from
$G$. But then $\hat H^*(G,k)$ has an invertible element of nonzero
degree, which implies that $G$ has rank one.}.  It follows that the
nil radical is not finitely generated, so $\hat H^*(G,k)$ is not
noetherian.
\end{proof} 

When $\hat H^*(G,k)$ is noetherian, modulo its radical it is a graded field, so the theory of supports coming from this ring is not interesting. This is why we chose for $R$ the noetherian subring $H^*(G,k)$
of $\hat H^{*}(G,k)$.  It is shown in \cite{BK:2002a} that contraction of ideals gives a one-one correspondence between prime ideals $\fp$ in $\hat H^*(G,k)$ and prime ideals $\fp^c$ in $H^*(G,k)$ in the nonperiodic case, and that the injective hulls $E(\hat H^*(G,k)/\fp)$ and $E(H^*(G,k)/\fp^c)$ are equal when $\fp$ is not maximal.

For $\sfT=\StMod kG$ and $R=H^*(G,k)$, the theory developed in Section
\ref{se:tensor} of this paper coincides with the theory developed by
Benson, Carlson and Rickard in \cite{BCR:1996a}. We should like to
note that the maximal ideal $\fm = H^+(G,k)$ of $R$ is not in the
support of any object; this follows, for example, from
\eqref{eq:recol}. However, $\fm$ is in the support of the cohomology
of some modules, such as the trivial module $k$.  Moreover, the
supports of $H^*(G,M)$ and $\hat H^*(G,M)$ agree except possibly for
$\fm$, because $\hat H^-(G,M)$ is $\fm$-torsion.

The correspondence of notation between \cite{BCR:1996a} and this
article is as follows. We denote by $V_G$ the maximal ideal spectrum
of $H^*(G,k)$ which is a homogeneous affine variety. Let $V$ be an
irreducible subvariety of $V_G$ corresponding to a homogeneous prime
ideal $\fp$ of $H^*(G,k)$. Let $\mcW$ be a specialization closed set
of homogeneous prime ideals in $H^*(G,k)$ and identify $\mcW$ with the
set of closed homogeneous subvarieties $W\subseteq V_G$ whose
irreducible components correspond to ideals in $\mcW$.

\begin{center}
\begin{tabular}{cc}
  BCR \cite{BCR:1996a} & This paper \\ \hline
  $k$ &  $\bbi$ \\
  $H^*(G,k)$ & $R$ \\
  $E(\mcW)$ & $\gam_{\mcW}\bbi$ \\
  $F(\mcW)$ & $L_{\mcW}\bbi$ \\
  $E(V)$ & $\gam_{\mcV(\fp)}\bbi$ \\
  $F'(V)$ & $L_{\mcZ(\fp)}\bbi$ \\
  $\kappa_V=E(V) \otimes_k F'(V)$ & $\gam_\fp\bbi$ \\
  $\kappa_V \otimes_k M$ & $\gam_\fp M$ \\ 
  $\mcV_G(M)$ & $\supp_R M$ \\ \hline \\
\end{tabular}
\end{center}

The kappa module $\kappa_V$ of \cite{BCR:1996a} is denoted $\kappa_\fp$ in some other papers. The crucial isomorphisms $E(\mcW)\cong\gam_{\mcW}\bbi$ and $F(\mcW)\cong L_{\mcW}\bbi$ follow from Theorem~\ref{thm:localization-finite}. We keep our correspondence of notation, and restate Theorem~\ref{thm:gammap-variants} to obtain the following. 

\begin{theorem}
\pushQED{\qed}%
Let $\mcV$ and $\mcW$ be specialization closed sets of homogeneous prime ideals in $H^*(G,k)$ such that $\mcV\setminus\mcW=\{\fp\}$.  Then
\[  
  E(\mcV) \otimes_k F(\mcW) \cong \kappa_V\,. \qedhere
\]
 \end{theorem} 

One feature of the support variety theory for $\StMod kG$ that does not hold more generally is the tensor product theorem: the support of the tensor product of two modules is the intersection of their supports.  This is Theorem~10.8 of \cite{BCR:1996a}. If, instead of using the whole cohomology ring, we just use a subring, then the tensor product theorem fails, even if the cohomology ring is a finitely generated module over the
chosen subring. To see this, just choose a subring such that there exist two distinct prime ideals in $H^*(G,k)$ lying over the same prime in the subring, and tensor the kappa modules for these primes. The tensor product is projective, but the intersection is not empty.

Conjecture 10.7.1 of \cite{Benson:2004b} was an attempt to find a way to compute the support variety $\mcV_G(M)$ from the cohomology of $M$.  This conjecture is false; see Example \ref{eg:V8} below. However, the following theorem does compute the support variety from the cohomology, and may be thought of as a replacement for this conjecture. 

\begin{theorem}
\label{thm:Bensonconj}
Let $\fp$ be a nonmaximal homogeneous prime ideal in $H^*(G,k)$ and let $\zeta_i\in H^{n_i}(G,k)$ ($1\le i\le s$) be nonzero elements such that $\fp=\sqrt{(\zeta_1,\dots,\zeta_s)}$. 
Let $L_{\zeta_i}$ be the kernel of a cocycle $\hat\zeta_i\colon \Omega^{n_i}k\to k$ representing $\zeta_i$.

Then for any $kG$-module $M$, the following are equivalent.
\begin{enumerate}[\quad\rm (1)]
\item $\fp \in \mcV_G(M)$.
\item There exists a simple $kG$-module $S$ such that
\[ 
\Ext^*_{kG}(S\otimes_k L_{\zeta_1} \otimes_k \dots \otimes_k L_{\zeta_s},M)_\fp\ne 0\,. 
\]
\end{enumerate}
\end{theorem}

\begin{proof}
This follows from Theorem~\ref{thm:detectionplus}. The Koszul object
$\kos S{(\zeta_1,\dots,\zeta_s)}$ in this context is, up to a shift,
the tensor product $S \otimes_k L_{\zeta_1}\otimes_k \dots\otimes_k
L_{\zeta_s}$.
\end{proof}  

Note that, in the preceding result, we could have also used Tate cohomology modules
$\widehat \Ext^*_{kG}(-,-)$. 
 
\subsection*{The derived category}
Another choice for $\sfT$ is  $\sfD(\Mod kG)$, the derived category. This is a rather poor choice, for the
following reason. If $G$ is a $p$-group, then the only localizing subcategories of $\sfD(\Mod kG)$ are zero and the whole category. Even if $G$ is not a $p$-group, the structure of the set of localizing subcategories in no way reflects the set of prime ideals in $H^*(G,k)$. The other problem in this case is that the trivial
module $k$, regarded as an object in $\sfD(\Mod kG)$, is not compact. So although there is a tensor product and there are function objects, the unit for the tensor product is not compact, and we cannot apply
the theory of Section~\ref{se:tensor}. 

\subsection*{The homotopy category of complexes of injective modules}
The third choice for $\sfT$ is the homotopy category of complexes of injective $kG$-modules, denoted $\sfK(\Inj kG)$.  This tensor triangulated category is investigated in Benson and Krause \cite{BK:2006}, and in particular there is a recollement
\begin{equation}
\label{eq:recol}
\StMod kG \simeq \sfK_\ac(\Inj kG) \ 
\begin{smallmatrix} \funct_k(tk,-) \\
\hbox to 50pt{\leftarrowfill} \\ \hbox to 50pt{\rightarrowfill} \\ 
\hbox to 50pt{\leftarrowfill} \\ - \otimes_k tk
\end{smallmatrix} \ 
\sfK(\Inj kG) \ \begin{smallmatrix} \funct_k(pk,-) \\
\hbox to 50pt{\leftarrowfill} \\ \hbox to 50pt{\rightarrowfill} \\ 
\hbox to 50pt{\leftarrowfill} \\ - \otimes_k pk
\end{smallmatrix} \ 
\sfD(\Mod kG). 
\end{equation}

Here, $ik$, $pk$ and $tk$ denote an injective, projective, Tate resolution of $k$ as a $kG$-module respectively, and $\sfK_\ac(\Inj kG)$ denotes the full subcategory of $\sfK(\Inj kG)$ consisting of
acyclic complexes. The equivalence between this and the stable module category is given by the theory of Tate resolutions of modules. The recollement follows from Theorem~\ref{thm:recollement} because $\sfK_\ac(\Inj kG)$ consists of all objects $X$ in $\sfK(\Inj kG)$ such that $\supp_RX$ does not contain the maximal ideal $\fm$. 

The compact objects in $\sfK(\Inj kG)$ are the semi-injective resolutions, see \cite{AFH}, of finite complexes of finitely generated $kG$-modules, and the full subcategory of compact objects forms a triangulated category equivalent to $\sfD^b(\mod kG)$. The object $ik$ is compact, and is the unit for the tensor product, so this plays the role of $\bbi$.  The function objects of Section \ref{se:tensor} are provided by Hom complexes
\[ 
\funct(X,Y)_n = \prod_p\Hom_k(X_p,Y_{n+p}) 
\] 
with the usual differential $(d(f))(x)=d(f(x))-(-1)^{|f|}f(d(x))$ and $G$-action $(g(f))(x)=g(f(g^{-1}x))$. If $C \to X$ and $D \to Y$ are semi-injective resolutions of objects $C$ and $D$ in $\sfD^b(\mod kG)$ then the map $\funct(C,D) \to \funct(C,Y)$ is a semi-injective resolution, and the map $\funct(X,Y) \to \funct(C,Y)$ is an isomorphism in $\sfK(\Inj kG)$.  It follows that for compact objects $X$ and $Y$, the function complex
$\funct(X,Y)$ is again compact. In particular, $X^\vee=\funct(X,ik)$ is a semi-injective resolution of $\funct(C,k)$.  It follows that compact objects are strongly dualizable.  The graded endomorphism ring of $ik$ is $H^*(G,k)$, which is our choice for $R$. The theory of varieties for modules in this context agrees with the theory set up in section 9 of \cite{BK:2006}.  It restricts to the full subcategory $\StMod kG$ of acyclic complexes to give the theory of \cite{BCR:1996a}. Exactly one more prime ideal comes into play, namely
the maximal ideal $\fm$, and this reflects the right hand side of the recollement. 

\begin{example}
\label{eg:V4}
We imitate the example in Section~\ref{se:commutative} to give an
example in $\sfK(\Inj kG)$ where the triangulated support differs from
the cohomological support. This gives a rather natural looking
example. We then modify it in Example \ref{eg:V8} to provide a
counterexample to Conjecture 10.7.1 of Benson~\cite{Benson:2004b}. Let
$G=(\bbz/2)^2$ and $k$ a field of characteristic two.  We have
$R=H^*(G,k)=k[x,y]$, where the generators $x$, and $y$ have degree
one. Let $Q$ be the (homogeneous) field of fractions of $R$. Then the
minimal injective resolution of $R$ has the form
\begin{equation}
\label{eq:KInj-resol} 
\cdots \to 0 \to Q \to \bigoplus_{\height\fp=1} E(R/\fp) \to E(R/\fm) \to 0 \to \cdots  
\end{equation}
where $\fm=H^+(G,k)$ is the unique maximal (homogeneous) ideal in $R$.  Recall from \cite[\S10]{BK:2006} that there is a functor $T$ from injective $R$-modules to $\sfK(\Inj kG)$ such that $H^*(G,T(I)) \cong I$. Apply this functor to the first nontrivial arrow in the resolution \eqref{eq:KInj-resol} and complete to a triangle in $\sfK(\Inj kG)$ to obtain
\begin{equation}
\label{eq:X}
 X \to T(Q) \to \bigoplus_{\height\fp=1}T(E(R/\fp))\to. 
\end{equation}
Then take cohomology to deduce that there is a short exact sequence
\[ 
0 \to E(R/\fm)[-1] \to H^*(G,X) \to R \to 0. 
\]
This sequence splits---for instance, the right hand side is a free module---to give
\begin{equation}
\label{eq:HX}
H^*(G,X) \cong R \oplus E(R/\fm)[-1].
\end{equation}
In fact, it turns out that $X$ is nothing other than $tk$, a Tate resolution for $k$ as a $kG$-module, and $H^*(G,X) \cong \hat H^*(G,k)$.  Using \eqref{eq:X} and \eqref{eq:HX}, we see that $\fm$ is in $\supp_R H^*(G,X)$ but not in $\supp_R X$.
\end{example} 

\begin{example}
\label{eg:V8}
We modify Example \ref{eg:V4} to give a counterexample, in $\StMod kG$, to Conjecture 10.7.1 of Benson \cite{Benson:2004b} as follows. Let $G=(\bbz/2)^3$ and let $k$ be a field of characteristic two. Then $R=H^*(G,k)=k[x,y,z]$, where the generators $x$, $y$ and $z$ have degree one. Let $\mcZ=\{\fq\in\spec R\mid\fq\not\subseteq(x,y)\}$.  Working in $\sfK(\Inj kG)$, we write $F=L_\mcZ(k)$, so that $H^*(G,F)=R_{(x,y)}$ is a homogeneous local ring with maximal ideal $(x,y)$. The minimal injective resolution of $R_{(x,y)}$ over $R$ has the form
\[ 
\cdots \to 0 \to Q \to 
\bigoplus_{\substack{\height\fp=1 \\ \fp\subseteq(x,y)}}
E(R/\fp) \to E(R/(x,y)) \to 0 \to \cdots 
\]
where $Q$ is the (homogeneous) field of fractions of $R$. Apply $T$ to the first nonzero arrow and complete to a triangle in $\sfK(\Inj kG)$ to obtain
\begin{equation}
\label{eq:X2} 
X \to T(Q) \to \bigoplus_{\substack{\height\fp=1 \\ \fp\subseteq(x,y)}} T(E(R/\fp))\to. 
\end{equation}
Applying cohomology, we deduce that there is an exact sequence of $R_{(x,y)}$-modules
\[ 
0 \to E(R/(x,y))[-1] \to H^*(G,X) \to R_{(x,y)} \to 0,. 
\]
This splits to give an isomorphism
\begin{equation}
\label{eq:HX2}
H^*(G,X) \cong R_{(x,y)} \oplus E(R/(x,y))[-1]. 
\end{equation}
Using \eqref{eq:X2}, we see that $\supp_R X = \{(0)\} \cup \{\fp\mid\height\fp=1,\fp\subseteq(x,y)\}$ while
using \eqref{eq:HX2}, we see that $\supp_R H^*(G,X)=\{\fp\mid\fp\subseteq(x,y)\}$. Now the maximal ideal $\fm$ of $H^*(G,k)$ is not in $\supp_R X$, so we may regard $X$ as an object in $\StMod(kG)$ using the recollement 
\ref{eq:recol}. Namely, $X$ is an acyclic complex, and the kernel of the middle differential gives an object $M$ in $\StMod(kG)$ whose triangulated support and cohomological support differ. 

\begin{theorem}
For $G = (\bbz/2)^3$ and $k$ a field of characteristic two, there exists a module $M$ in $\StMod(kG)$ such that $\mcV_G(M)=\{(0)\}\cup \{\fp\mid\height\fp=1,\fp\subseteq(x,y)\}$, so that $(x,y)\not\in \mcV_G(M)$, but $(x,y)\in \supp_{H^*(G,k)} H^*(G,M)$. \qed
\end{theorem}
\end{example} 

Subsequently, families of such examples have been constructed in \cite{BC:2007}.

\subsection*{Finite Dimensional Algebras.}
Generalizing the theory for finite groups, there is a theory of support varieties for modules over a finite dimensional algebra $A$, developed by Snashall and Solberg~\cite{SS}; see also Solberg~\cite{So}.  The
theory is developed there for finitely generated $ A$-modules and for the derived category $\sfD^b(\mod A)$ of bounded complexes of finitely generated modules. Their idea is to use the natural homomorphism from the Hochschild cohomology ring $\HH^*( A)$ to $\Ext^*_{ A}(M,M)$.  In the case where $M$ is the quotient of $A$ by its radical, they show that the kernel of this map consists of nilpotent elements.  
A problematic feature of the theory is that it is not known in general whether the quotient of $\HH^*(A)$ by its nil radical is finitely generated as an algebra. This issue is discussed at length in \cite{SS,So}.  

The way to use our theory to construct support varieties for modules over finite dimensional algebras is to use the triangulated category $\sfT=\sfK(\Inj A)$, whose compact objects form a copy of $\sfD^b(\mod A)$. There is a natural homomorphism of rings $\HH^*(A) \to Z(\sfT)$.  If $R$ is a finitely generated subalgebra of $\HH^*(A)$ then we may apply our theory to the composite $R \to \HH^*( A) \to Z(\sfT)$.  This way, we extend the Snashall--Solberg theory to objects in $\sfK(\Inj A)$.  

If $A$ is self-injective the theory is further refined in \cite{EHSST} using the stable category of finitely generated modules $\stmod  A$.  The comments in the group case about the relationship between the large stable module category $\StMod A$ and $\sfK(\Inj A)$ apply equally well here: one has homomorphisms $R \to \HH^*( A) \to Z(\StMod A)$ and a theory of varieties for infinitely generated modules in this situation. So for example our Theorem~\ref{thm:KRS} generalizes Theorem 7.3 of \cite{EHSST} to infinitely generated modules. 

For a finite group scheme $G$ over an arbitrary field $k$, Friedlander and Pevtsova developed in \cite{FP} the notion of support using so-called $\pi$-points. They define for each $kG$-module $M$ a support space $\Pi(G)_M$, which can be identified with a set of nonmaximal homogeneous prime ideals of the cohomology ring $H^*(G,k)$.  This approach generalizes the work from \cite{BCR:1996a} for finite groups. In particular, we have $\Pi(G)_M=\supp_RM$ as before when we take $\sfT=\StMod kG$ and $R=H^*(G,k)$.

\section{Complete intersections}
\label{se:ci}
In this section we discuss support for complexes of modules over commutative complete intersection rings. The main goal is to show only how the theory presented here relates to the one of Avramov and Buchweitz~\cite{Av:vpd,AB} for finitely generated modules. Further elaborations are deferred to another occasion. 

Let $A$ be a commutative noetherian ring of the form $Q/I$, where $Q$
is a commutative noetherian ring and $I$ is an ideal generated by a
regular sequence; see \cite[Section~16]{Mat} for the notion of a regular
sequence.  The principal examples are the local complete intersection
rings; see Remark~\ref{rem:AB-theory}.  In what follows $\sfK(A)$
denotes the homotopy category of complexes of $A$-modules.

Let $R$ be the ring of cohomology operators defined by the surjection $Q\to A$, see~\cite{AB:two,AS}, and the paragraph below. Thus, $R=A[\chi_1,\dots,\chi_c]$, the polynomial ring over $A$ in variables $\chi_i$ of degree $2$. For each complex $X$ of $A$-modules, there is a natural homomorphism of rings $R\to \Ext^*_A(X,X)$, with image in the center of the target. These define a homomorphism from $R$ to the center of the derived category of $A$. It lifts to an action on $\sfK(A)$ and gives a homomorphism of rings $R\to \cent(\sfK(A))$.  This map is the starting point of everything that follows in this section, so we sketch a construction. It is based on a method in \cite{AB:two} for defining cohomology operators, and uses basic notions and constructions from Differential Graded homological algebra, for which we refer to \cite[\S2]{Av:barca}.  

\subsection*{Cohomology operators}
Let $Q\to A$ be a homomorphism of commutative rings.  Let $B$ be a semi-free resolution of the $Q$-algebra $A$. Thus, $B$ is a DG algebra over $Q$ whose underlying algebra is the graded-symmetric algebra over a graded free $Q$-module, concentrated in cohomological degrees $\leq0$, and there is a quasi-isomorphism $\eps\col B\to A$ of DG $Q$-algebras. Set $B^{\sfe}=B\otimes_QB$; this is also a DG $Q$-algebra, and, since $B$ is
graded-commutative, the product map $\mu\col B^{\sfe}\to B$, where $b'\otimes b''\mapsto b'b''$, is a morphism of DG algebras.  Let $U$ be a semi-free resolution of $B$ viewed as a DG module over $B^{\sfe}$ via $\mu$. Set $V= A\otimes_B U$; this has the structure of a DG module over $A$, that is to say, a complex of $A$-modules. 
Let $X$ be a complex of $A$-modules; it acquires a structure of a DG $B$-module via $\eps$. Associativity of tensor products yields isomorphisms of complexes of $A$-modules
\[  
  X\otimes_B U \cong X\otimes_A(A\otimes_BU) = X\otimes_A V\,.
\]
The functors $X\otimes_B -$ and $-\otimes_A V$, on the categories of DG modules over $B^{\sfe}$ and over $A$, respectively, are additive and hence preserve homotopies.  One has homomorphisms of algebras
\[  
  \Hom_{\sfK(B^{\sfe})}^*(U,U) \xra{\ X\otimes_B -\ }
  \Hom_{\sfK(A)}^*(X\otimes_A V,X\otimes_A V)  \xleftarrow{\ -\otimes_A V\ }
  \Hom_{\sfK(A)}^*(X,X)\,.
\]
The complex of $A$-modules $V$ is semi-free and the map $V= A\otimes_BU\to A\otimes_BB=A$ is a quasi-isomorphism, and hence it is a homotopy equivalence. Therefore, the map $-\otimes_AV$ in the diagram above is an isomorphism. Consequently, one obtains the homomorphism of algebras
\[  
  \Hom_{\sfK(B^{\sfe})}^*(U,U)\lto   \Hom_{\sfK(A)}^*(X,X)\,.
\]
It is clear that this map is natural in $X$, so it induces a homomorphism of rings
\[  
  \Hom_{\sfK(B^{\sfe})}^*(U,U) \lto  Z(\sfK(A))\,.
\]
We note that $\Hom_{\sfK(B^{\sfe})}^*(U,U)$ is the Shukla
cohomology~\cite{Sh} 
of $A$ over $Q$; it is also called the Hochschild-Quillen cohomology.

When $A=Q/(\bsq)$, where $\bsq=q_1,\dots,q_c$ is a $Q$-regular sequence, a standard calculation, see \cite[(2.9)]{AB:two}, yields an isomorphism of $A$-algebras
\[  
  \Hom_{\sfK(B^{\sfe})}^*(U,U)\cong A[\chi_1,\dots,\chi_c]\,,
\]
where the $\chi_i$ are indeterminates in cohomological degree $2$. They are the cohomology operators defined by the presentation $Q\to A$. This completes our sketch. 

\subsection*{Support}
Let $A=Q/(\bsq)$, where $\bsq=q_1,\dots,q_c$ is a $Q$-regular sequence and  $R=A[\chi]$ the ring of cohomology operators, and $R\to \cent(\sfK(A))$ the homomorphism of rings, introduced above.  Set $\sfK=\sfK(\Inj A)$, the homotopy category of complexes of injective $A$-modules; it is a compactly generated triangulated category, see~\cite[Prop.~2.3]{Kr:stable}.  Restricting $R\to \cent(\sfK(A))$ gives a homomorphism of rings
\[  
  A[\chi]=R \lto \cent(\sfK)\,.
\]
One has thus all the ingredients required to define local cohomology and support for objects in $\sfK$. 

\begin{remark}
\label{rem:kvsd}  
Via the standard embedding of the derived category $\sfD(A)$ of $A$
into $\sfK$, this then applies to all complexes of $A$-modules.  To be
precise: Restriction of the quotient functor $\sfq\col\sfK(A)\to
\sfD(A)$ admits a fully faithful right adjoint, say $\sfi$; it maps $X$ to 
a semi-injective resolution of $X$. These functors fit into the following diagram of functors:
\[  
\xymatrixcolsep{2pc} 
\xymatrixrowsep{2.5pc} 
\xymatrix{ 
   \sfK^c \ar@{->}[d]^\inc \ar@<-1ex>[rrr]_-{\sfq}^-{\sim} &&& 
     \sfD^f(A) \ar@{->}[d]_\inc \ar@<-1ex>[lll]_-{\sfi}\\
      \sfK
       \ar@<-1ex>[rrr]_-{\sfq}
    &&& \sfD(A) \ar@<-1ex>[lll]_-{\sfi}}
\]
The equivalence in the top row holds by \cite[Prop.~2.3]{Kr:stable}. 

Henceforth, when we talk about the support of a complex $X$ in $\sfD(A)$, we mean $\supp_R\sfi X$.  In the diagram, $\sfD^f(A)$ denotes the subcategory $\sfD(A)$ consisting of complexes with finitely generated cohomology.  In particular, one arrives at a notion of support for complexes in $\sfD^f(A)$.  In Remark~\ref{rem:AB-theory}, we compare this construction with the one of Avramov and Buchweitz. First, we record an observation.
\end{remark} 

For a homomorphism $A\to R$ of commutative rings, the \emph{fibre} at
a point $\fp$ in $\spec A$ is the ring $R\otimes_Ak(\fp)$. The lemma
below says that the support of any finitely generated $R$-module is
detected along its fibers.

\begin{lemma}
\label{lem:fibres}
Let $A\to R$ be a homomorphism of commutative noetherian rings. For each finitely generated $R$-module $L$, one has an equality
\[  
  \supp_RL = \bigsqcup_{\fp\in\spec A}\supp_{(R\otimes_Ak(\fp))}(L\otimes_Ak(\fp))\,.
\]
\end{lemma} 

\begin{proof}
The crucial remark is that for each finitely generated module $M$ over a commutative noetherian ring $B$, one has $\supp_BM=\{\fq\in\spec B\mid M_\fq\ne0 \}$; see Lemma~\ref{le:supp-ann}. Now fix a prime ideal $\fp$ in $ A$. Since $L$ is finitely generated as an $R$-module, $L\otimes_Ak(\fp)$ is finitely generated as an $R\otimes_Ak(\fp)$-module. This remark, and the isomorphism $L\otimes_A k(\fp)\cong L\otimes_R(R\otimes_Ak(\fp))$, yield
\[  
  \supp_{(R\otimes_Ak(\fp))}(L\otimes_Ak(\fp)) = \supp_RL \cap \supp_R(R \otimes_Ak(\fp))\,.
\]
The upshot is that it suffices to prove the desired equality for $L=R$, and in this case it is evident.
\end{proof} 

We now present one of the main results in this section.  Recall that a complex of modules over a ring is said to be \emph{perfect} if it is quasi-isomorphic (in the derived category) to a finite complex of finitely generated projective modules. 

\begin{theorem}
\label{thm:ci-fibres}
Let $A=Q/(\bsq)$, where $\bsq=q_1,\dots,q_c$ is a $Q$-regular sequence, and set $\sfK=\sfK(\Inj A)$. Let $R$ be the induced ring of cohomology operators on $\sfK$. 

For each $X$ in $\sfK^{c}$ which is perfect over $Q$, one has $\supp_RX = \supp_R\End^*_{\sfK}(X)$. Furthermore, there is a fibre-wise decomposition
\[  
  \supp_RX = \bigsqcup_{\fp\in\spec A}\supp_{(R\otimes_Ak(\fp))}\Ext^*_{A_\fp}(X_\fp,k(\fp))\,.
\]
\end{theorem} 

\begin{remark}
When a complex $X$ in $\sfK^{c}$ has  finite injective dimension over $Q$, then again one has $\supp_RX=\supp_R\End^*_{\sfK}(X)$ and
\[  
  \supp_RX = \bigsqcup_{\fp\in\spec A}\supp_{(R\otimes_Ak(\fp))}\Ext^*_{A_\fp}(k(\fp),X_\fp)\,.
\]
The proof is similar to the one for Theorem~\ref{thm:ci-fibres}.
\end{remark} 

\begin{proof}[Proof of Theorem~\emph{\ref{thm:ci-fibres}}]
When $Y$ is a compact object in $\sfK$, it is semi-injective, so one has an identification $\Hom_{\sfK}^{*}(-,Y)=\Ext_{A}^{*}(-,Y)$. A crucial result in the proof is that, for such a $Y$, the $R$-module $\Ext_A^*(X,Y)$ is finitely generated. This holds because $X$ is perfect over $Q$ and the $A$-module $H^*(Y)$ is finitely generated; see Gulliksen~\cite[(2.3)]{Gu}, and also \cite[(5.1)]{AS}, \cite[(4.2)]{AGP}. This fact is used implicitly, and often, in the argument below.  For each prime $\fp$ in $\spec A$, we write $\sfK_{\fp}$ for $\sfK(\Inj A_\fp)$.  
Theorem~\ref{thm:finiteness}(2) yields the first equality below:
\begin{align*}
\supp_R X &= \supp_R\End^*_{\sfK}(X) \\  
     &= \bigsqcup_{\fp\in\spec A}\supp_{(R\otimes_Ak(\fp))}(\End^*_{\sfK}(X)\otimes_Ak(\fp))\\
     &= \bigsqcup_{\fp\in\spec A}\supp_{(R\otimes_Ak(\fp))}
           (\End^*_{\sfK_{\fp}}(X_\fp)\otimes_{A_\fp}k(\fp))\,.
\end{align*}
The second one holds by Lemma~\ref{lem:fibres}, while the third holds because the $A$-module $H^*(X)$ is finitely generated. The homomorphism of rings $R\to Z(K)$ yields a homomorphism of rings $R_\fp\to Z(K_\fp)$.  Note that $R_\fp$ is a ring of cohomology operators defined by the presentation $A_\fp =Q_{\wt \fp}/(\bsq)_{\wt\fp}$, where $\wt\fp$ is the preimage of $\fp$ in $Q$. Furthermore, $X_\fp$ is perfect over $Q_{\wt \fp}$. 

Therefore we may assume $A$ and $Q$ are local rings, with residue field $k$.  The desired result is then that for each complex $X$ of $A$-modules which is perfect over $Q$, one has an equality
\[  
  \supp_R(\Ext^*_A(X,X)\otimes_Ak)  = \supp_R \Ext^*_A(X,k)\,.
\]
The set on the left contains the one on the right, since the ring $R$ acts on $\Ext^*_A(X,k)$ through $\Ext^*_A(X,X)$.  We prove that the equality on the left holds. 

In the rest of the proof, it is convenient to set $X(-)=\Ext_A^*(X,-)$, viewed as a functor from $\sfD^f(A)$, the derived category of cohomologically finite complexes of $A$-modules, to the category of finitely generated graded $R$-modules.  We repeatedly use the characterization of support in Lemma~\ref{le:supp-ann}(1). 

Let $\fm$ be the maximal ideal of $A$, and consider the full subcategory
\[  
  \sfC = \{M\in \sfD^f(A)\mid \mcV(\fm)\cap \supp_R X(M) \subseteq \supp_R X(k)\}
\]
As $\supp_R(X(M)\otimes_Ak) = \mcV(\fm R) \cap \supp_R X(M)$, the desired result follows from: 

\begin{Claim}
$\sfC=\sfD^f(A)$
\end{Claim} 

Indeed, since $\supp_{R}X(k)$ is closed $\sfC$ is a thick subcategory of $\sfD^f(A)$; this is easy to check.  It thus suffices to prove that $\sfC$ contains all finitely generated $A$-modules. This is now verified by a standard induction argument on Krull dimension. To begin with, observe that $k$ is in $\sfC$, and hence so is any module of dimension zero, for such a module admits a finite filtration with subquotients isomorphic to $k$. Let $M$ be a finitely generated $A$-module with $\dim M\geq 1$. Let $L$ be the $\fm$-torsion submodule of $M$, and consider the exact sequence of finitely generated $A$-modules $0\to L\to M\to M/L\to 0$. Note that $L$ is in $\sfC$, since $\dim L=0$, so to prove that $M$ is in $\sfC$ it suffices to prove that $M/L$ is in $\sfC$. Thus, replacing
$M$ by $M/L$ one may assume that there is an element $a\in \fm$ which is a nonzero divisor on $M$. Consider the exact sequence
\[  
  0\lto M \xra{\ a\ } M \lto M/aM\lto 0\,.
\]
This gives rise to an exact sequence of $R$-modules $X(M)\xra{a} X(M)\to X(M/aM)$, so one deduces that 
\[  
  \mcV(aR)\cap \supp_R X(M) = \supp_R (X(M)/aX(M)) \subseteq \supp_RX(M/aM)\,.
\]
Since $\dim(M/aM) = \dim M-1$ the induction hypothesis yields that $M/aM$ is in $\sfC$. Given this, the inclusion above implies $M$ is in $\sfC$. 

This completes the proof of the claim, and hence of the theorem.
\end{proof} 

\begin{remark}
\label{rem:AB-theory}
Let $A$ be a complete intersection local ring, with maximal ideal $\fm$ and residue field $k$. Assume $A$ is $\fm$-adically complete. Cohen's Structure Theorem provides a presentation $A\cong Q/I$ with $(Q,\fq,k)$ a regular local ring and $I\subseteq \fq^2$.  Since $A$ is a complete intersection, the ideal $I$ is generated by a regular sequence of length $c=\dim Q - \dim A$; see \cite[\S21]{Mat}.  Let $A[\chi]=A[\chi_1,\dots,\chi_c]$ be the corresponding ring of operators.  Let $M$ be a finitely generated $A$-module. Observe that any cohomologically finite complex of $A$-modules is perfect over $Q$, since $Q$ has finite global dimension.  Theorem~\ref{thm:ci-fibres} and Lemma~\ref{le:supp-ann}(1) show that the fibre of $\supp_R\sfi M$, see Remark \ref{rem:kvsd}, at the maximal ideal $\fm$ is precisely the subset
\[  
  \supp_{k[\chi]}\Ext^*_A(M,k) = \mcV(\fa) \subseteq \spec k[\chi]\,,
\]
where $\fa=\ann_{k[\chi]}\Ext^*_A(M,k)$. Let $\wt k$ be the algebraic closure of $k$. The Nullstellensatz implies that the sets $\mcV(\fa)$ and 
\[  
  \{(b_1,\dots,b_c)\in \wt k^c\mid f(b_1,\dots,b_c)=0\text{ for } f\in \fa\}\cup \{0\}
\]
determine each other. The latter is precisely the support variety of $M$ in the sense of Avramov~\cite{Av:vpd}; see also Avramov and Buchweitz~\cite{AB}.
\end{remark} 

\begin{remark}
\label{rem:AB-compare}  
Let $A=Q/I$ where $Q$ is a commutative noetherian ring of finite global dimension and $I$ is generated by a regular sequence. Let $A[\chi]$ be the associated cohomology operators. Theorem~\ref{thm:ci-fibres} applies to each compact object in $\sfK$.  It follows from this result and Remark~\ref{rem:AB-theory} that for each cohomologically finite complex $X$ of $A$-modules one has a natural projection $\supp_RX\to \spec A$, and the fibre over each prime $\fp$ in $\spec A$ encodes the support of the complex of $A_\fp$-modules $X_\fp$, as defined in \cite{AB}. Observe that this set is empty outside $\supp_AX$, in the sense of Section~\ref{se:commutative}.  Our notion of support is thus a refinement of the one in \cite{Av:vpd} by means of the classical support. 

The discrepancy between the two is clarified by considering the perfect complexes over a local ring $A$, with residue field $k$. The support in the sense of \cite{AB} of any nonzero perfect complex $P$ is then $\{0\}$, which corresponds to the ideal $(\chi)\subset k[\chi]$. On the other hand, $\supp_RP$ is a closed subset of $\spec A=\mcV(\chi)\subset \spec A[\chi]$ and equals the classical support $\supp_AP$.  This information is important if one wants to classify the thick subcategories of $\sfD^f(A)$. Bear in mind that thick subcategories of perfect complexes are classified by specialization closed subsets of $\spec A$, see~\cite{Ne:dr}. 

There are other noteworthy aspects to the construction in this section, the most important one being that it gives a local cohomology theory, with supports in the ring of cohomology operators, for complexes over complete intersections. Here is one evidence of its utility: Specializing Corollary~\ref{cor:connectedness} yields a connectedness theorem for support varieties of complexes over complete intersection rings. Using this, one can recover a result of Bergh~\cite[(3.2)]{Be}; the details will be provided elsewhere.
\end{remark} 

\begin{ack}
During the course of this work, the second author was collaborating with Lucho Avramov on a related project~\cite{AI}. We would like to thank him, and also Jesse Burke, Amnon Neeman, and a referee, for criticisms and suggestions regarding this write-up. The second and third authors are grateful to the Mathematisches Forschungsinstitut at Oberwolfach for vital support through a `Research in Pairs' stay. The first author is grateful to the Humboldt Foundation for enabling extended visits to the third author.
\end{ack} 

\bibliographystyle{amsplain}

\begin{thebibliography}{99} 

\bibitem{JVS1}
{\sc L. Alonso Tarrío, A. Jerem\'ias L\'opez, M. J. Souto Salorio:} Localization in categories of complexes and unbounded resolutions, Canad. J. Math.  \textbf{52}  (2000),  225--247. 

\bibitem{JVS2}
{\sc L. Alonso Tarrío, A. Jerem\'ias L\'opez, M. J. Souto Salorio:} Bousfield localization on formal schemes. J. Algebra \textbf{278} (2004),  585--610. 

\bibitem{Av:vpd}{\sc L. L. Avramov:} Modules of finite virtual projective dimension, Invent. Math.  \textbf{96}  (1989), 71--101.  

\bibitem{Av:barca}{\sc L. L. Avramov:} Infinite free resolutions.  Six lectures on commutative algebra (Bellaterra, 1996), 1--118, Progr. Math., \textbf{166}, Birkh\"auser, Basel, 1998.   

\bibitem{AB:two}{\sc L. L. Avramov and R.-O. Buchweitz:} Homological algebra modulo a regular sequence with special attention to codimension two. J. Algebra \textbf{230} (2000), 24--67. 

\bibitem{AB}{\sc L. L. Avramov and R.-O. Buchweitz:} Support varieties and cohomology over complete intersections.  Invent. Math.  \textbf{142} (2000),  285--318. 

\bibitem{AFH}{\sc L. L. Avramov, H.-B. Foxby and S. Halperin:} Differential graded homological algebra, in preparation.  

\bibitem{AGP}{\sc L. L. Avramov, V. N. Gasharov, and I. V. Peeva:} Complete intersection dimensions, Inst. Hautes \'Etudes Sci. Publ. Math.  \textbf{86} (1997), 67--114. 

\bibitem{AS} {\sc L. L. Avramov and  L.-C. Sun:} Cohomology operators defined by a deformation,
J. Algebra \textbf{204} (1998), 684--710. 

\bibitem{AI} {\sc L. L. Avramov and S. B. Iyengar:} Modules with prescribed cohomological support, Ill. Jour. Math. \textbf{51} (2007), 1--20. 

\bibitem{Ba}{\sc P. Balmer:} Supports and filtrations in algebraic geometry and modular representation theory, Amer. J. Math. \textbf{129} (2007), 1227--1250.

\bibitem{BBD}{\sc A. Beilinson, J. Bernstein, and P. Deligne:} Faisceaux pervers,   Ast\'erisque \textbf{100}, Soc. Math. France, 1983. 

\bibitem{Benson:1991}{\sc D. J. Benson:} Representations and cohomology II. Cambridge studies in advanced mathematics, \textbf{31}, Cambridge Univ. Press, Cambridge, 1991. 

\bibitem{Benson:2004b}{\sc D. J. Benson:} Commutative algebra in the cohomology of groups. Trends in commutative algebra, 1--50, Math. Sci. Res. Inst. Publ., \textbf{51}, Cambridge Univ. Press, Cambridge, 2004. 

\bibitem{BC:2007} {\sc D. J. Benson and J. F. Carlson:} Varieties and cohomology of infinitely generated modules, preprint 2007.

\bibitem{BCR:1996a} {\sc D. J. Benson, J. F. Carlson, and J. Rickard:} Complexity and varieties for infinitely generated modules. II.  Math. Proc. Cambridge Philos. Soc. \textbf{120} (1996), 597--615. 

\bibitem{BIK:2008} {\sc D. J. Benson, S. B. Iyengar, and H. Krause:} Localizing subcategories of the stable module category of a finite group, in preparation.

\bibitem{BK:2002a} {\sc D. J. Benson and H. Krause:} Pure injectives and the spectrum of the cohomology ring of a finite group.  J. Reine Angew. Math.  \textbf{542} (2002), 23--51. 

\bibitem{BK:2006} {\sc D. J. Benson and H. Krause:} Complexes of
injective $kG$-modules, Algebra Number Theory, to appear.

\bibitem{Be}{\sc P. Bergh:} On support varieties of modules over complete intersections, Proc. Amer. Math. Soc. \textbf{135} (2007), 3795--3803.

\bibitem{BH}{\sc W. Bruns and J. Herzog:} Cohen-Macaulay rings.  Cambridge Studies in Advanced Mathematics, \textbf{39}. Cambridge University Press, Cambridge, 1998.  Revised edition. 

\bibitem{Ca:1983}{\sc J. F. Carlson:} The varieties and the cohomology ring of a module. J. Algebra \textbf{85} (1983), 104--143. 

\bibitem{Ca}{\sc J. F. Carlson:} The variety of an indecomposable module is connected. Invent. Math.  \textbf{77} (1984),  291--299. 

\bibitem{CE}{\sc H. Cartan and S. Eilenberg:} Homological algebra. Princeton Univ. Press, Princeton, NJ, 1956. 

\bibitem{Ch}{\sc S. Chebolu:} Krull-Schmidt decompositions for thick subcategories. J. Pure Appl. Algebra
\textbf{210}  (2007), 11--27.

\bibitem{EHSST}{\sc K. Erdmann, M. Holloway, N. Snashall, \O. Solberg and R. Taillefer:} Support varieties for selfinjective algebras.  $K$-Theory \textbf{33} (2004), 67--87. 

\bibitem{Fo}{\sc H.-B. Foxby:} Bounded complexes of flat modules. J. Pure Appl. Algebra \textbf{15} (1979), 149--172. 

\bibitem{FI}{\sc H.-B. Foxby and S. Iyengar:} Depth and amplitude for unbounded complexes. Commutative Algebra (Grenoble/Lyon 2001), Contemp. Math. \textbf{331}, A.M.S., R.I., 2003; 119--137. 

\bibitem{FP}{\sc E. Friedlander and J. Pevtsova:} $\Pi$-supports for modules for finite group schemes over a field. Duke Math. J. \textbf{139} (2007), 317--368.

\bibitem{GZ}{\sc P.\ Gabriel and M.\ Zisman:} Calculus of fractions and
homotopy theory.  Ergebnisse der Mathematik und ihrer Grenzgebiete
{\bf 35}, Springer-Verlag, New York (1967).

\bibitem{GM}
{\sc J. P. C. Greenlees and J. P. May:}  Derived functors of $I$-adic completion and local homology.  J. Algebra \textbf{149}  (1992), 438--453. 


\bibitem{Gr}{\sc A. Grothendieck:} Cohomologie locale des faisceaux coh\'erents et th\'eor\`emes de Lefschetz locaux et globaux (SGA 2).  S\'eminaire de G\'eom\'etrie Alg\'ebrique du Bois Marie, 1962.  Revised reprint of the 1968 French original. Documents Math\'ematiques (Paris) \textbf{4}. Soc. Math. France, Paris, 2005. 

\bibitem{Gu}{\sc T. H. Gulliksen:} A change of ring theorem with applications to Poincar\'e series and intersection multiplicity.  Math. Scand.  \textbf{34} (1974), 167--183. 

\bibitem{Ha} {\sc R. Hartshorne:} Residues and duality. Lecture Notes Math. \textbf{20}, Springer-Verlag, New York (1966). 

\bibitem{Ho} {\sc M.Hopkins:} Global methods in homotopy theory, in: Homotopy theory (Durham, 1985), London Math. Soc. Lecture Note Ser. \textbf{117}, Cambridge Univ. Press, Cambridge, 1987, 73--96.
 
\bibitem{HPS} {\sc M. Hovey, J. Palmieri, and N. Strickland:} Axiomatic stable homotopy theory.  Mem. Amer. Math. Soc.  \textbf{610}, Amer. Math. Soc., (1997). 

\bibitem{IK}{\sc S. Iyengar and H. Krause:} Acyclicity versus total acyclicity for complexes over noetherian rings. Documenta Math. \textbf{11} (2006), 207--240. 

\bibitem{Kr:krs}{\sc H. Krause:} Decomposing thick subcategories of the stable module category.  Math. Ann.  \textbf{313} (1999), 95--108. 

\bibitem{Kr:br}{\sc H. Krause:} A Brown representability theorem via coherent functors.  Topology \textbf{41} (2002), 853--861. 

\bibitem{Kr:stable}{\sc H. Krause:} The stable derived category of a noetherian scheme. Compos. Math.  \textbf{141} (2005), 1128--1162. 

\bibitem{Kr:thick}{\sc H. Krause:} Thick subcategories of modules over
commutativer rings, Math. Ann., to appear.

\bibitem{LMS}{\sc L. G. Lewis Jr., J. P. May, M. Steinberger, and J. E. McClure:} Equivariant stable homotopy theory. Lecture Notes Math. \textbf{1213}, Springer-Verlag, Berlin, 1986.

\bibitem{Li}{\sc J. Lipman:} Lectures on local cohomology and duality.  Local cohomology and its applications (Guanajuato, Mexico), Lecture Notes. Pure Appl. Math. \textbf{226}, Marcel Dekker, New York  (2002), 39--89. 

\bibitem{MacL}{\sc S. MacLane:} Categories for the working mathematician. Springer Verlag (1971).
 
\bibitem{Ma}{\sc H. R. Margolis:} Spectra and the Steenrod Algebra. North-Holland (1983). 

\bibitem{Mat}{\sc H. Matsumura:} Commutative ring theory. Cambridge University Press (1986). 

\bibitem{Ne}{\sc A. Neeman:} Triangulated categories, Annals of Math. Studies \textbf{148}, Princeton Univ. Press, Princeton, NJ, 2001. 

\bibitem{Ne:dr}{\sc A. Neeman:} The chromatic tower of $D(R)$. Topology \textbf{31} (1992), 519--532.  

\bibitem{Qu}{\sc D. Quillen:} The spectrum of an equivariant cohomology ring. I, II, Ann. of Math. \textbf{94} (1971), 549--572; ibid. \textbf{94} (1971), 573--602.  

\bibitem{Ri}{\sc J. Rickard:} Idempotent modules in the stable category. J. London Math. Soc. \textbf{178} (1997), 149--170. 

\bibitem{Sh}{\sc U. Shukla:} Cohomologie des alg\`ebres
associatives. Ann. Sci. \'Ecole Norm. Sup. \textbf{78} (1961)
163--209.

\bibitem{SS} {\sc N. Snashall and \O. Solberg:} Support varieties and Hochschild cohomology rings,
Proc. London Math. Soc. (3) \textbf{88} (2004), 705--732. 

\bibitem{So} {\sc \O. Solberg:} Support varieties for modules and complexes, Trends in representation theory of algebras and related topics, Contemp. Math. \textbf{406}, Amer. Math. Soc., Providence, RI, 2006; 239--270. 

\bibitem{Th}
{\sc R. Thomason:} The classification of triangulated subcategories, Compositio Math. \textbf{105} (1997), 1--27.  

\end{thebibliography}

\end{document}